%% file: main.tex
\documentclass[a4paper]{article}

\usepackage[english]{babel}
\usepackage[utf8]{inputenc}
\usepackage{amsfonts}
\usepackage{booktabs}
\usepackage{tabu}
\usepackage[T1]{fontenc}

\usepackage[a4paper,top=1cm,bottom=1.5cm,left=2cm,right=2cm,marginparwidth=1cm]{geometry}

\usepackage{amsmath}
\usepackage{graphicx}
\usepackage{caption}
\usepackage[colorinlistoftodos]{todonotes}
\usepackage[colorlinks=true, allcolors=blue, breaklinks]{hyperref}
\usepackage{wrapfig}
\usepackage{dblfloatfix}
\usepackage{titling}
\usepackage{float}

\usepackage{algorithmic}
\usepackage[ruled,vlined]{algorithm2e}
\usepackage{fancyvrb}

\usepackage{url}

\usepackage{breakurl}

\usepackage{tikz}
\usetikzlibrary{shapes.geometric, arrows}

\begin{document}

\title{\Large \textbf{A Simulation Framework for Ride-Hailing with Electric Vehicles}}

\author{
    \normalsize \textbf{Chen Zhang$^1$},
    \normalsize \textbf{Sushil Mahavir Varma$^2$}\\
    \\
    \begin{minipage}{0.9\textwidth}
    \centering
    \small $^1$ Georgia Institute of Technology \\
    \small $^2$ INRIA Paris 
    \end{minipage}
}

\date{}
\maketitle

\vspace{-0.5cm}
\begin{abstract}

    This research presents a Python-based simulation framework designed to model electric vehicle (EV) on-demand transportation systems, with a focus on optimizing urban fleet operations. Built on a process-driven architecture, the system efficiently simulates EV fleet dynamics, including passenger matching, vehicle dispatching, and charging strategies, while enabling customization to address critical challenges such as charger placement, fleet management, and algorithm performance. We overcome the challenge of high dimensional state-space and non-Markovian system dynamics by executing processes asynchronously using SimPy and updating only the states that are affected by employing object-oriented programming. As a result, our simulation framework is capable of handling peak demand scenarios involving thousands of trips and completing multi-day scenarios in minutes. The modular design enables users to experiment with parameters, test algorithms, and integrate custom datasets, making the tool highly adaptable for diverse urban contexts. By providing a realistic and extensible platform, this adaptable, scalable, and open-source framework advances the optimization of EV fleet operations and offers a valuable resource for decision-makers and city planners navigating the transition to sustainable urban mobility solutions. 
    
    We also present a case study using the NYC taxi dataset evaluating various dispatching algorithms, including closest vehicle dispatch, closest available vehicle dispatch, and power-of-d vehicle dispatch, and exploring charging approaches like continuous and nighttime charging. We propose a novel \emph{adaptive} power-of-$d$ dispatch policy, which dynamically adjusts to real-time conditions and demonstrates high throughputs when combined with adaptive charging policies that interrupt charging to meet demand during the peak and delay some of the charging to the nighttime.
    
\end{abstract}

\vspace{0.5cm}

\input{paper}

\bibliographystyle{unsrt}
\bibliography{refs}

\end{document}

%% file: paper.tex

\section{Introduction}


\noindent The transportation industry is in the midst of a shift from Internal Combustion Engine cars to Electric Vehicles due to technology advancements\cite{ford2021}, growing climate change awareness, and widespread governmental backing\cite{padilla2020}. This transformation is manifesting in the on-demand transportation industry, where ride-hailing companies provide electric vehicle (EV) services and plan to switch to a completely EV-based fleet\cite{alto2023, uber2020, lyft2023, bellan2021}. In turn, it becomes essential for an EV-based on-demand transportation system to have a comprehensive understanding of the electric vehicle design, the infrastructure required\cite{lamonaca2022state}, and design efficient matching and charging algorithms to deliver a good quality of service.

\noindent We aim to simulate a fleet of autonomous electric vehicles with a spatial queueing mechanism for an urban transportation network. Our simulation involves a certain amount of EVs and chargers placed throughout the city. The customer demand for rides can be modeled based on different datasets, such as NYC taxi data\cite{nyc2023} and Chicago data\cite{chicago2023}. Once a ride request from a customer is received, the matching algorithm in the system evaluates and decides which EV to dispatch and serve the customer. If no EV is appropriate to meet the trip, the customer's request cannot be fulfilled and will be dropped. Once an EV is dispatched, it travels to pick up the customer and then proceeds to the customer's destination. Throughout their operations to serve the trips, EVs lose their State of Charge (SoC) and require frequent charging. In turn, the charging algorithm determines the most effective timing and location for each EV to recharge in order to sustain operations. This initiative is motivated by the advent of electric vehicle fleets deployed as taxis. Companies like Blusmart\cite{business2024}, Alto\cite{dallas2024}, and Revel\cite{boylan2023} have been pioneers in integrating electric vehicles into their fleet operations.

\noindent The code for the simulation framework described in this paper is publicly available on GitHub \cite{varma_zhang_sim}.

\subsection{Operational Challenges in Electric Vehicle Fleet Management}

\noindent Integrating electric vehicle fleets poses a myriad of operational challenges in the landscape of urban mobility. There exists a substantial amount of analytical work, but those models are generally simplified and have a gap with reality. City planners and fleet operators need to make effective decisions to manage the EV fleets to make the operations efficient. 
There are various questions related to the operations of an electric vehicle fleet that need to be addressed.

\noindent \textbf{Fleet Parameters:} How to set the fleet parameters, such as the fleet size and the battery pack capacity? The right fleet size is vital to balance service levels and costs. Too few vehicles can lead to long wait times, failing to meet customer demand. Conversely, an excessively large fleet can increase operational costs due to unused capacity. Similarly, setting the appropriate battery pack capacity is key to optimizing service efficiency. Large batteries offer a longer travel range but are more expensive, while small batteries are cheaper but need more frequent recharging sessions because of limited travel range.

\noindent \textbf{Charger Parameters:} How many charging stations are needed? How many charging posts should each station include? What should be the operating charge rate for these posts? These decisions directly impact the efficiency of charging EVs. Too many chargers and posts will increase the costs, whereas too few chargers limit charging options, forcing vehicles to travel further to find available stations. Insufficient charging posts fail to meet the demand, causing long wait times. In addition, the charge rate affects how quickly vehicles can be charged and returned to service. A high charge rate enables quick charging but requires expensive infrastructure, while a low charge rate is cheaper but increases charging time.

\noindent \textbf{Matching Algorithms:} The matching algorithms determine how to pair an EV with an incoming customer request. These algorithms not only have to identify which vehicles are available within a given time frame and region but also select the most suitable vehicle from those available ones to serve the customer arrival. Well-designed matching algorithms optimize the efficiency of the ride-hailing service by pairing customers with the most suitable EVs to minimize wait times, maximize vehicle usage, and diminish operational costs.

\noindent \textbf{Charging Algorithms:} The charging algorithms are established to specify the optimal battery level at which an EV should be directed to recharge, determine the best times to send vehicles to chargers, and set the duration that each car should spend charging. Well-built charging algorithms ensure that vehicles are consistently powered for maximum availability, thereby preventing service delays and enhancing the overall efficiency of the fleet. The timing of charging can also be adjusted based on the demand intensity. For example, charging can be increased during periods of low demand to ensure that EVs are fully prepared to meet higher demand when it arises.

\noindent \textbf{Charging Stations Locations:} Where should the chargers be located in the city? This is a complex problem because the optimal charger placement should consider factors like traffic patterns, population density, proximity to high-demand locations, etc. If the chargers are not put strategically, it can result in underutilized resources, increased congestion in certain areas, and inefficient charging for electric vehicles. Specifically, this can cause some EVs to wait a long time to be charged, while other chargers remain unused for much of the time.

\noindent \textbf{Car Initial Locations:} How to distribute the fleet at the start of each day? To initialize the car locations to meet the demand, random placement might not be effective because different regions have different levels of demand. For instance, residential neighborhoods may require more EVs in the morning as people commute to work, whereas commercial districts may have higher demand in the evening as people return home.

\subsection{Main Contributions}

\noindent To this end, we propose the creation of a simulation architecture for EV-based transportation systems that closely mimics real-world dynamics. This open-source platform considers and integrates all critical questions outlined above, providing decision-makers with a tool for modeling and analyzing the performance of EV fleets under various matching and charging algorithms. This tool empowers planners to experiment with different scenarios and enhance the efficiency of EV fleet operations. This manuscript is a companion to \cite{varma_ev} and describes the simulator used and developed in the process.

\noindent The simulator operates on a process-based discrete-event simulation framework using SimPy, where each discrete event activates some specific processes. Each customer arrival initiates a sequence of events: an EV picks up the customer, and other EVs may be directed to charging stations. These events then trigger further actions in the simulation. For example, once the pickup process is completed, it then initiates the process of driving with customers to their destination, followed by additional related processes. 

\noindent Our simulation framework offers the following advantages.

\subsubsection{Curse of Dimensionality}

\noindent \textbf{High-dimensional state space of the system:} Our system is designed to monitor the activities of all EVs, so we keep track of such a huge and high-dimensional state space. It contains the locations of all EVs, the SoCs of all EVs, the states of all EVs (driving with passenger, picking up passenger, charging, waiting at the charger, driving to the charger, idle), and the origin and destination of all active trips.

\noindent However, we do not continuously monitor the locations and SoCs of all EVs. Instead, we use a ``teleportation'' approach where these states are updated only once a process is finished. For instance, we do not keep track of the location constantly while an EV is serving a trip but only update it after the trip is completed. Additionally, there exist lots of attributes for all EVs. It is impractical to update all of them each time when a process is finished. Therefore, we selectively update only those attributes that are relevant to the completed process by employing Object-Oriented Programming in Python. In particular, each EV is an object of a class with attributes like current state, current SoC, and current location. Then, these attributes are asynchronously updated across the fleet. In particular, at the end of each process in the simulation, we update the attributes of the EVs that are affected by that process. Such an asynchronous update is time-efficient as opposed to updating the state of the entire system after every process.

\subsubsection{Non-Markovianity of the System}

\noindent The system dynamics are non-Markovian with respect to the state space considered. To make it Markovian, we would need to continuously keep track of the locations and SoCs of all EVs at all times. At any given time in the simulation, several different processes are running in parallel. For example, while some of the EVs are serving requested trips, others may be driving to the charger or actively charging. Each of these processes runs for a stipulated time and may also trigger other processes in the meantime. For instance, an EV that is currently charging may be matched to a customer request, which then leads to an interruption of the charging process and triggers the process of the EV driving to the customer's pickup location. Implementing multiple such processes in an asynchronous fashion in a time-efficient manner poses a challenge.

\noindent We addressed this obstacle by using the SimPy library in Python, which is well-suited for handling non-Markovian dynamics. Particularly, in SimPy, we can add timeouts for a stipulated amount of time which can be used to simulate events. For example, to model the process of picking up the customer, we schedule a timeout corresponding to the time that the EV takes to reach the passenger. When the timeout concludes, we then update the location and SoC. Another advantage of SimPy is that such timeouts can be triggered for different objects in the simulation in parallel. In this way, we can run several processes in parallel, such as EVs picking up passengers while others are charging or driving to the charger. Lastly, we can also use this framework to trigger timeouts as a result of other timeouts. This can be used to trigger processes as a result of the completion of other processes. 

\subsubsection{Modularity}
\noindent Our modular simulation architecture offers flexibility, allowing for comprehensive analysis for all questions mentioned above. The practitioners have the freedom to adjust any parameters and test the performance of various algorithms. Users can modify settings including chargers locations, initial EV locations, fleet characteristics and charger behaviors. The framework supports the users to develop and test custom matching and charging algorithms tailored to meet specific operational needs. They also have the option to import their own dataset into the simulation, enabling insights that closely mirror real-world conditions. We provide the detailed instructions on how users can incorporate custom algorithms and datasets in the GitHub README file.

\noindent We conduct a comprehensive study by varying different parameters of the simulation model and prescribing the rules of the questions above. Our study also includes tests of multiple matching and charging algorithms that have been proposed in the literature.

\subsubsection{Scalability}
\noindent Our system is capable of efficiently handling peak demand scenarios involving 3,000 to 10,000 trips, with a corresponding number of electric vehicles. For the simulations presented in later sections, the system completes each one-day simulation in approximately 10 minutes (600-650 seconds) and each three-day simulation in around 40 minutes (2200-2500 seconds). Users also have the flexibility to run simulations for any number of days while maintaining fast runtimes.

\section{Literature Review}

\noindent Existing literature explores various aspects of ride-hailing platforms, including capacity planning, and control algorithm design for pricing, matching, relocation, and charging. Specifically, research has shown different approaches to capacity planning related to queueing systems \cite{halfin1981heavy} and spatial models \cite{besbes2022spatial}. The design of near-optimal algorithms for pricing, matching, and relocation in the context of ride-hailing has been proposed in the literature \cite{banerjee2015pricing,banerjee2018state,mahavir2020dynamic,varma2021dynamic,ozkan2020dynamic}. Furthermore, studies on managing large-scale charging stations \cite{wang2016smart,shao2023preemptive} and scheduling charging of electric vehicles \cite{wu2022smart} have enriched our understanding of operational strategies in electric vehicle-based ride-hailing systems. These analyses are generally conducted by dividing the space into discrete grids \cite{bimpikis2019spatial, braverman2019empty} and analyzing the resultant system.

\noindent However, the existing models often oversimplify the complex dynamics of real-world scenarios. Moreover, there is no unified simulation framework that can be easily adapted to study varios algorithms from the literature. Recognizing these limitations, our objective is to develop a flexible, real-world simulation platform. This platform is implemented to design and test matching, charging, and pricing algorithms for EV-based ride-hailing systems, bridging the gap between theoretical research and practical application.

\noindent There exist several papers that focus primarily on autonomous vehicles (AVs) without accounting for the dynamics of electric vehicles. For instance, Levin et al. (2017) present a model that incorporates dynamic ride-sharing into shared autonomous vehicle (SAV) models by integrating their behavior into traffic flow simulations to highlight its efficiency in reducing congestion and travel times, without considering aspects of EV operation. Vosooghi et al. (2019b) utilize multi-agent simulations to demonstrate the influence of fleet size, vehicle capacity, and rebalancing strategies on SAV service efficiency. The subsequent study by Vosooghi et al. (2019a) incorporates user trust and willingness to use SAVs, affecting fleet sizing and consequently demand patterns, service efficiency, and operational metrics. These studies delve into the broader aspects of AV systems, rather than electric vehicle management and operations.

\noindent Bauer et al. (2018) employ an agent-based model to analyze optimal fleet battery ranges, necessary charging infrastructure, and fleet sizes for shared automated electric vehicles (SAEVs) operating as taxis, emphasizing cost savings and environmental benefits. Loeb and Kockelman (2019) evaluate the performance and cost-effectiveness of SAEV fleets compared to hybrid-electric vehicles, highlighting environmental benefits of SAEVs despite lower profitability. Yang et al. (2023) study the design and evaluation of electric autonomous mobility-on-demand systems, considering urban congestion and limited link space with a discrete event simulation model, aiming to optimize fleet size and charging station locations to mitigate congestion effects and enhance vehicle utilization. These studies align with our focus on studying fleet operations through simulation models, but they restrict their analysis to fixed matching and charging algorithms and concentrate on specific objectives such as fleet sizing and charging infrastructure.

\noindent By contrast, our system subsumes a wide range of algorithms that can be tested and compared. Our simulation framework offers extensive customization options and allows for the testing of various matching and charging algorithms, providing users with the tools to explore a broad range of operational parameters. For example, our framework supports interruptible charging processes and flexible charging schedules, which can lead to more efficient fleet utilization that are not included in the cited studies. The flexibility makes our simulation framework valuable for researchers and planners seeking to test and evaluate EV fleet operations in varied urban environments.

\section{Overview of the Simulation Framework}

Our simulation framework is based on SimPy, the discrete event simulator library which allows us to perform non-Markovian simulations efficiently and also incorporate various real-life features. In particular, we track the state of charge (SoC) of all electric vehicles in addition to their locations. We also take into account the pick-up times, travel time to the charging stations, and queueing at the charger.

\subsection{State of the System:}

\noindent For each EV, the system keeps track of its location, current state, and State of Charge (SoC). There are six possible states for an EV: (1) Driving with Passenger: The EV is actively transporting a passenger to their destination. (2) Picking up Passenger: The EV is en route to pick up a passenger. (3) Charging: The EV is plugged in and charging at a station. (4) Waiting at the Charger: The EV is queued at a charging station but not yet charging. (5) Driving to the Charger: The EV is on its way to a charging station to recharge. (6) Idle: The EV is not currently engaged in any of the above activities.

\noindent Instead of continuously monitoring the EV's location as it moves, the system only updates the car's location once it reaches the destination. Similarly, this method is also employed for monitoring the SoC, where updates are only made after the charging process is completed. In this way, the need for constant monitoring and data transmission is reduced.

\noindent The system records the number of EVs occupying the charging stations. There are two primary states for a charger: (1) Available: The charger has at least one charging post that is free. (2) Busy: All posts of the charger are occupied, indicating full capacity. By default, a charging station is considered busy if the number of cars at the charging station is greater than or equal to the number of posts of that charging station. Otherwise, the station is available.

\noindent The trip requests of the simulation go through several states, and this is a matchmaking process between passengers and EVs. There are four states for a trip: (1) Unavailable: No EV is available to fulfill the trip request. (2) Waiting: The trip request has been made, and the system is searching for an available EV. (3) Matched: A suitable EV has been successfully assigned to the trip arrival, and the pickup process will begin or has already started. (4) Reneged: The trip request could not be fulfilled or matched.



\subsection{Physics of the System}



\begin{figure}[H]
\centering
\usetikzlibrary{decorations.markings}
\begin{tikzpicture}[node distance=1cm, auto,
    cross/.style={
        decoration={
            markings,
            mark=at position 0.17cm with with {\fill[black] circle (3pt);},
            mark=at position 0.67cm with {\fill[black] circle (3pt);},
            mark=at position 1.67cm with {\fill[black] circle (3pt);},
            mark=at position 2.57cm with {\fill[black] circle (3pt);},
            mark=at position 3.27cm with {\fill[black] circle (3pt);},
            mark=at position 4.47cm with {\fill[black] circle (3pt);},
            mark=at position 5.17cm with {\fill[black] circle (3pt);},
            mark=at position 5.97cm with {\fill[black] circle (3pt);},
            mark=at position 7.17cm with {\fill[black] circle (3pt);},
        },
        postaction={decorate}
    }]

  \tikzstyle{process} = [rectangle, minimum width=2cm, minimum height=1cm, text centered, draw=black, fill=orange!30, align=center]
  \tikzstyle{arrow} = [thick,->,>=stealth, align=center]

  \node (pro1) [process] {Charging\\Algorithm};
  \node (pro2) [process, right=of pro1] {Driving to\\Charger};
  \node (pro3) [process, right=of pro2, xshift=1cm] {Queueing at\\Charger};
  \node (pro4) [process, right=of pro3, xshift=1cm] {Charging};
  \node (pro5) [process, below=of pro1, yshift=-1.5cm] {Matching\\Algorithm};
  \node (pro6) [process, above right=of pro5, xshift=1.5cm, yshift=-1cm] {Interrupt};
  \node (pro7) [process, below right=of pro5, xshift=5cm, yshift=1cm] {Pick Up\\Passenger};
  \node (pro8) [process, right=of pro7] {Drop Off\\Passenger};
  \node (start) [left=of pro5] {};

  \draw[->, cross] (-2.3,-6) -- (12,-6) node[right] {t};
  \draw [arrow] (start) |- (pro1);
  \draw [arrow] (start) |- (pro5);
  \draw (start) -- ++(0,-2.5);
  \draw (start) -- ++(0,-1.5) -- ++(0.5,0) -- ++(0,-1);
  \draw (start) -- ++(0,-1.5) -- ++(1.5,0) -- ++(0,-1);
  \draw (start) -- ++(0,-1.5) -- ++(2.4,0) -- ++(0,-1);
  \draw (start) -- ++(0,-1.5) -- ++(3.1,0) -- ++(0,-1);
  \draw (start) -- ++(0,-1.5) -- ++(4.3,0) -- ++(0,-1);
  \draw (start) -- ++(0,-1.5) -- ++(5,0) -- ++(0,-1) -- ++ (1, -1);
  \draw (start) -- ++(0,-1.5) -- ++(5.8,0) -- ++(0,-1);
  \draw (start) -- ++(0,-1.5) -- ++(7,0) -- ++(0,-1) node[midway, right] {...};
  \draw [arrow] (pro1) -- (pro2);
  \draw [arrow] (pro2) -- (pro3) node[midway, above, align=center] {If not\\terminated};
  \draw [arrow] (pro3) -- (pro4) node[midway, above, align=center] {If not\\terminated};
  \draw [arrow] (pro5) -- ++(2.5,0) node[midway, above] {Car Idle?} |- node[near end, above] {no} (pro6);
  \draw [arrow] (pro5) -- ++(2.5,0) |- node[near end, above] {yes} (pro7);
  \draw [arrow] (pro6) -| (pro7) node[near start, below, align=center] {Update SoC,\\State, Location};
  \draw [arrow] (pro7) -- (pro8);
  \draw [arrow] (pro6.north) -- ++(0,1) -| (pro2.south) node[near start, below left] {If driving to charger,\\terminate};
  \draw [arrow] (pro6.north) -- ++(0,1) -| (pro3.south) node[near start, below right] {If queueing at charger,\\terminate};
  \draw [arrow] (pro6.north) -- ++(0,1) -| (pro4.south) node[near end, below right] {If charging,\\terminate};
  \draw [arrow] (pro4.east) -- ++(2,0) node[midway, above] {Update\\SoC, State};
  \draw [arrow] (pro8.east) -- ++(2.5,0) node[midway, above] {Update SoC,\\State, Location};

\end{tikzpicture}

\vspace{-0.8cm}
\begin{figure}[H]
    \centering
    \hspace{2cm}
    \includegraphics[width=6cm]{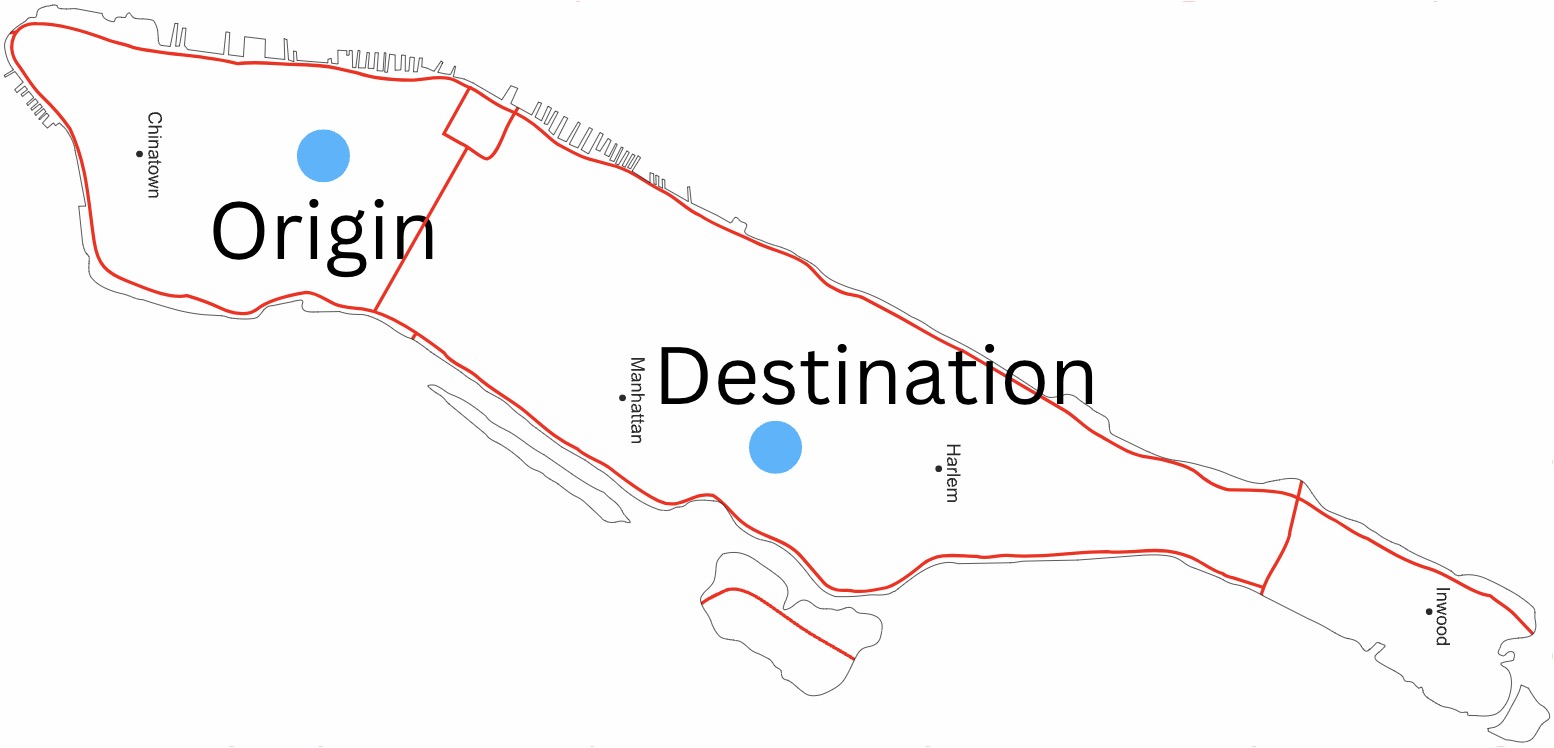}
    \label{fig:transition_diagram}
\end{figure}

\vspace{-3.5cm}
\noindent \hspace{-10cm} customer arrivals
\vspace{2.5cm}
\caption{Flow Chart of the System}
\label{fig:arrivals_matching_charging}
\end{figure}

\noindent As shown in the Fig.\ref{fig:arrivals_matching_charging}, the flow chart illustrates the primary logic of the simulation. The time axis displayed at the bottom tracks the arrival times of customers. The customer arrivals either follow a Poisson process or are based on a predefined dataset. Each customer arrival is associated with a specific origin and destination within the city boundary. Upon a customer's arrival, both matching and charging algorithms are triggered.

\noindent The matching algorithm assigns a vehicle to serve the customer. If the selected vehicle is not idle—by default meaning it is en route to a charger, waiting in line at a charger, or currently charging—the simulation interrupts the vehicle's current activity, updates its SoC and location, and sets its state to idle. If the vehicle is already idle, it picks up the passenger and transports them to the designated destination. Once the passenger is dropped off, the car's SoC and location are refreshed, and it returns to an idle state.

\noindent Concurrently, the charging algorithm sends a list of EVs to chargers. These vehicles will then drive to the charger, queue at the charger, and proceed to charge. If a vehicle is called upon to serve a customer during this process, its current activity will be terminated. If not, the car will continue until it finishes charging. Afterward, the SoC and the state of the EV are updated.

\subsubsection{Matching}
\begin{figure}[H]
\centering
\begin{tikzpicture}[node distance=1cm, auto]
  \tikzstyle{process} = [rectangle, minimum width=2cm, minimum height=1cm, text centered, draw=black, fill=orange!30, align=center]
  \tikzstyle{arrow} = [thick,->,>=stealth, align=center]
  \node (pro1) [process] {Matching\\Algorithm};
  \node (start) [left=of pro1, xshift=-4cm] {};
  \node (end) [right=of pro1, xshift=2cm] {};
  \draw [arrow] (start) |- (pro1) node[near end, above, align=center] {Location of Customer,\\Location, SoC, State of EVs};
  \draw [arrow] (pro1) -- (end) node[midway, above, align=center] {EV to dispatch};
\end{tikzpicture}
\caption{Input and Output of Matching Algorithm}
\label{fig:matching_input_output}
\end{figure}
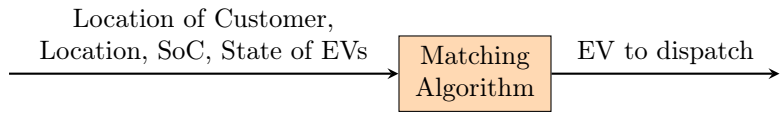

\noindent In the system, customer arrivals act as elementary events. As shown in the figure, for each arrival, the system makes a critical decision: whether to serve the customer or not. If the decision is to decline service, the customer promptly exits the system. If the system decides to serve the customer, a matching algorithm is used to identify the most suitable EV for the arrival. In particular, the matching algorithm takes into account the customer's location as well as EVs' locations and SoCs to determine which EV to dispatch for the pickup, as shown in Fig.\ref{fig:matching_input_output}.

\noindent Once the matching algorithm designates an EV to send to pick up the customer, the system first checks the state of the EV. If the EV is currently charging, waiting at the charger, or driving to the charger, its charging is interrupted by default so it can respond to the trip request. Once the EV is available to serve the customer, it navigates from its current location to the customer's pickup point, transports the customer to their destination, and finally completes the trip by dropping off the passenger. This process is illustrated in Fig.\ref{fig:arrivals_matching_charging}.

\subsubsection{Charging}
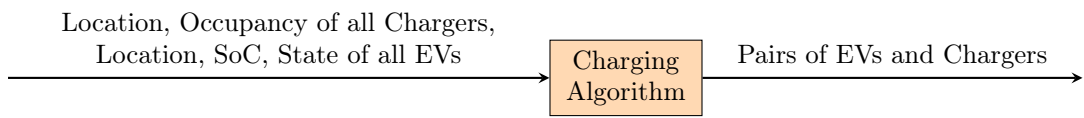
\begin{figure}[H]
\centering
\begin{tikzpicture}[node distance=1cm, auto]
  \tikzstyle{process} = [rectangle, minimum width=2cm, minimum height=1cm, text centered, draw=black, fill=orange!30, align=center]
  \tikzstyle{arrow} = [thick,->,>=stealth, align=center]
  \node (pro1) [process] {Charging\\Algorithm};
  \node (start) [left=of pro1, xshift=-6cm] {};
  \node (end) [right=of pro1, xshift=4cm] {};
  \draw [arrow] (start) |- (pro1) node[near end, above, align=center] {Location, Occupancy of all Chargers,\\Location, SoC, State of all EVs};
  \draw [arrow] (pro1) -- (end) node[midway, above, align=center] {Pairs of EVs and Chargers};
\end{tikzpicture}
\caption{Input and Output of Charging Algorithm}
\label{fig:charging_input_output}
\end{figure}

\noindent Meanwhile, when there is an arrival, the system also conducts real-time assessments of the EV fleet to determine the need for charging. This process is guided by the charging algorithm. The charging algorithm outputs two lists: a list of EVs that need to be sent to the chargers and a list of chargers that these EVs are assigned to. This is illustrated in Fig.\ref{fig:charging_input_output}.

\noindent Subsequently, each EV in the output list drives towards its assigned charger. Depending on the charging station's capacity and current occupancy, the EV either queues to wait for the next available post or starts charging immediately.

\noindent Throughout the charging phase, the EV remains connected to the charger. It keeps charging until either of two conditions is met: the EV's battery reaches full capacity, or the system dispatches the EV to go serving a customer. If the latter condition is triggered, the EV halts charging and promptly heads out to serve the customer. This is shown in Fig.\ref{fig:arrivals_matching_charging}.

\section{Details of the Simulation Model}

\subsection{Generating the Arrivals}

\noindent The dataset source for the simulation framework is categorized into two types: randomly generated data and real-life data.

\subsubsection{Randomly generated dataset}

\noindent This dataset is synthetically produced by the simulation framework itself. It is designed to mimic the real-world processes. The customer arrival times in the dataset are modeled using a Poisson process, which is used to describe the probability of a number of events happening within a fixed time interval. When the individual events occur randomly at a constant average rate, the time between each event follows an exponential distribution based on the arrival rate defined by the system. The locations, including the origins and destinations, are generated based on a uniform distribution across a spherical surface. This approach can ensure that the distances and locations are not biased towards any particular region.

\subsubsection{Real life dataset}

\noindent In contrast to the artificially generated data, this data source utilizes actual data that is collected from real-world operations. Our system does not impose any simplifying assumptions that might distort the complexity or variability observed in reality.

\noindent The system is designed with flexibility and applicability by allowing the users to import any dataset of their choice. The primary requirement is that the dataset should include specific information that is essential for simulation, including the start and end coordinates (latitude and longitude), trip distance, velocity, and time of each trip. This adaptability makes sure that the framework can be applied to a wide range of scenarios.

\noindent For demonstration, the system uses the New York Taxi dataset as an example. This dataset is chosen for its rich urban mobility data and realistic representation of taxi service operations.

\subsection{Processes associated with an EV}
The EV class includes the following processes associated with it.
\subsubsection{Pick up passenger}
The process of picking up a passenger begins by setting the trip state to \texttt{MATCHED}, which indicates that the customer has been successfully paired with the EV. Subsequently, the system calculates the distance from the current location to the pickup location and estimates the time required to reach the pickup point. Then, the car state is updated to \texttt{DRIVING\_WITHOUT\_PASSENGER}, and a timeout that equals to the estimated pickup time is set in the system. Upon reaching the pickup point, the car's location is updated to the trip's origin, which means the customer has been picked up, and the car state switches to \texttt{DRIVING\_WITH\_PASSENGER}.

\subsubsection{Drop off passenger}
When the vehicle reaches the pickup location, the system calculates the trip distance and total time to drive from the pickup location to the drop-off location. Concurrently, the vehicle begins its journey towards the passenger's destination. A corresponding timeout is set based on the calculated trip time. Upon reaching the destination, the car location is updated to match the trip's drop-off point. The energy consumption required to complete the trip is estimated by the system. In this way, the SoC of the vehicle's battery is adjusted downward to show the amount of energy spent while transporting the customer. Finally, the state of the car is changed to \texttt{IDLE}, which signals that the passenger is dropped off and the vehicle is ready for its next assignment.

\begin{algorithm}[H]
\caption{Run Trip} \label{alg:run_trip}
{\fontsize{10}{13}\selectfont
\SetAlgoLined
\SetKwInOut{Input}{Input}
\SetKwInOut{Initialize}{Initialization}
\Input{Trip}
\For{Every Arrival}{
\If{car state is \texttt{DRIVING\_TO\_CHARGER} or \texttt{CHARGING} or \texttt{WAITING\_FOR\_CHARGER}}{
\textnormal{Call Algorithm \ref{alg: interrupt_charging} to interrupt charging}\;
}
Change trip state to \texttt{MATCHED}\;
Calculate the pickup distance and trip pickup time\;
Change car state to \texttt{DRIVING\_WITHOUT\_PASSENGER}\;
Add a timeout equal to the trip pickup time\;
Calculate the trip time and trip distance\;
Set the car location equal to the trip's origin\;
Change car state to \texttt{DRIVING\_WITH\_PASSENGER}\;
Add a timeout equal to the trip time\;
Calculate the consumption to finish the trip and reduce the SOC accordingly\;
Set the car location equal to the trip's destination\;
Change car state to \texttt{IDLE}\;
}
}
\end{algorithm}

\subsubsection{Drive to charger}
When we send a car to the charger, the car's state is updated to \texttt{DRIVING\_TO\_THE\_CHARGER}. Simultaneously, a timeout is set for the estimated time duration of the drive from the vehicle's current location to its assigned charger. The system calculates the amount of energy consumption that is spent to drive to the charger, and the SoC of the car's battery decreases accordingly. Upon arrival at the charger, the EV's location is updated to be the charger's location, stating that it is now at the charging station. The car's state then changes to \texttt{WAITING\_AT\_THE\_CHARGER}. Following this, the Queueing at the Charger algorithm is invoked with the car's ID and the ending SOC as inputs.

\begin{algorithm}[H]
\caption{Drive to the Charger} \label{alg: drive_to_the_charger}
{\fontsize{10}{13}\selectfont
\SetAlgoLined
\SetKwInOut{Input}{Input}
\SetKwInOut{Initialize}{Initialization}
\Input{Charger ID, End SOC}
Change car state to \verb|DRIVING_TO_THE_CHARGER|\;
Add a timeout equal to the drive time to the charger\;
Reduce the SOC by the amount spent while driving to the charger\;
Set the car location equal to the charger location\;
Change car state to \verb|WAITING_AT_THE_CHARGER|\;
Call Algorithm \ref{alg: queueing} with inputs Car ID and End SOC\;
}
\end{algorithm}

\subsubsection{Charging}
After being removed from the queue of the charger, the EV is assigned to a charging post of the station and start charging. We assume that the charging posts are automated, meaning that cars start to charge as soon as they are allocated to a post. The charging process begins by changing the car's state to \texttt{CHARGING}. Meanwhile, a timeout is set corresponding to the anticipated time duration of the charging process. If the charging process is not interrupted, the vehicle's SoC is updated to match the end SoC, which is the battery's target charge level after the completion of charging. Once the car finishes charging, its state is reverted to \texttt{IDLE}, implying it is no longer charging and is ready to serve customers. At the same time, the occupancy of the charger is decremented by one, and the state of the charger is set to \texttt{AVAILABLE}.  

\begin{algorithm}[H]
\caption{Car Charging} \label{alg: charging}
{\fontsize{10}{13}\selectfont
\SetAlgoLined
\SetKwInOut{Input}{Input}
\SetKwInOut{Initialize}{Initialization}
\Input{Charger ID, End SOC}
Change car state to \verb|CHARGING|\;
Add a timeout equal to the charging time\;
Set the SOC equal to the End SOC\;
Change car state to \verb|IDLE|\;
Decrease the occupancy of the charger by one and set charger state to \verb|AVAILABLE|\;
Call Algorithm \ref{alg: queueing} so that cars waiting in the queue starts charging\;
}
\end{algorithm}

\subsubsection{Interrupt charging}
If the charging interruption functionality is enabled, based on the EV's current state, the system might interrupt the process when the EV is \texttt{DRIVING\_TO\_CHARGER}, \texttt{WAITING\_FOR\_CHARGER}, or \texttt{CHARGING}. If the vehicle is currently driving to the charger, the system first interrupts the process. The SoC and location are updated accordingly, and the car state is changed to \texttt{IDLE}. If the vehicle is waiting in the queue at a charging station, the system interrupts the process by changing the EV state to \texttt{IDLE} and removing the car from the queue list of that charger. If the vehicle is charging, the charging process is halted immediately. The EV's SoC is updated to show the increased battery level during the charging process before the interruption. Additionally, the occupancy of the charger is decremented by one after the interruption. The state of the charger is set to \texttt{AVAILABLE}, and the state of the EV is set to \texttt{IDLE}. After these updates, the Queueing at the Charger algorithm is triggered to allocate an available charging post to the next vehicle in line.

\begin{algorithm}[H]
\caption{Interrupt Charging} \label{alg: interrupt_charging}
{\fontsize{10}{13}\selectfont
\SetAlgoLined
\SetKwInOut{Input}{Input}
\SetKwInOut{Initialize}{Initialization}
\Input{Charger ID}
\uIf{car state is \texttt{DRIVING\_TO\_CHARGER}}
  {
    Interrupt the process corresponding to Algorithm \ref{alg: drive_to_the_charger}\;
    Update SOC, car location, car state to \texttt{IDLE};
  }
  \uElseIf{car state is \texttt{WAITING\_FOR\_CHARGER}}
  {
    Update car state to \texttt{IDLE}\;
    Remove the car from the charger queue list\;
  }
  \uElseIf{car state is \texttt{CHARGING}}
  {
    Interrupt the process corresponding to Algorithm \ref{alg: charging}\;
    Update SOC\;
    Decrease the occupancy by one\;
    Update car state to \texttt{IDLE}\;
    Update charger state to \texttt{AVAILABLE}\;
    Call Algorithm \ref{alg: queueing} so that cars waiting in the queue start charging\;
  }
}
\end{algorithm}

\subsection{Process Associated with a Charger}
The charger class includes the following processes associated with it.

\subsubsection{Queueing at the charger}
At the beginning of the process, the vehicle that will be charging at a charger is added to the queue of the charger, along with its target SoC. The process of managing the queue of cars waiting to charge keeps running if the charger state is not \texttt{busy} and the queue is not empty. If both conditions hold, the first car in the queue begins charging, and the system invokes the Car Charging algorithm. As a car starts charging, the number of occupied posts at the charger is incremented by one. The system then checks if this addition has occupied all available posts at the charging station. If so, the charger state is updated to \texttt{BUSY}, suggesting that no more cars can start charging until any post becomes available again. Meanwhile, the car that has just started charging is removed from the queue. This process repeats and continues to move the cars from the charger queue to the posts until either the charger becomes unavailable or there is no car waiting in the queue.

\begin{algorithm}[H]
\caption{Queueing at the Charger} \label{alg: queueing}
{\fontsize{10}{13}\selectfont
\SetAlgoLined
\SetKwInOut{Input}{Input}
\SetKwInOut{Initialize}{Initialization}
\Input{Car ID, End SOC}
Add input car and end SOC (if provided) to the charger queue\;
\While{Charger is not \texttt{BUSY} and charger queue is non-empty}{
First car of the queue starts charging (Call Algorithm \ref{alg: charging})\;
Increase the occupancy of the charger by one\;
\If{All posts are occupied}{
Set Charger state equal to \texttt{BUSY}\;
}
Remove the first car from the queue\;
}
}
\end{algorithm}


\section{Matching and Charging Algorithms} \label{sec: match_charge}

\subsection{Matching Algorithms}
\noindent We analyze the following three matching algorithms for dispatching electric vehicles that are adapted from \cite{varma_ev}: closest dispatch, closest available dispatch, and power-of-d vehicles dispatch. These algorithms focus on two key decisions: (1) Determining which EVs are available to be dispatched; (2) Choosing the most suitable EV to serve a customer request.

\subsubsection{First Decision: Identification of Available EVs}

\noindent The matching algorithm selects a subset of EVs and deems them to be available for dispatch. In our simulations, we experiment with the following four filtering criteria for defining available EVs:
\begin{itemize}
    \item Only Idle
    \item Idle, Charging, and Waiting at the Charger
    \item Idle, Charging, Waiting at the Charger, and Driving to the charger
    \item EVs that have been charging for a minimum stipulated time (e.g., for at least 10 minutes since they started charging)
\end{itemize}
Once we obtain the list of available EVs, the next step is to select one of these cars to dispatch to pick up the customer. This filtered list of available EVs is then used as input to the second decision stage.

\subsubsection{Second Decision: Selection of the EV for Dispatch}

\noindent To choose an EV to dispatch, we fetch the SoC and the location of the list of available EVs from the first decision. The second decision is made by either selecting an EV that can fulfill the customer's request or dropping the customer if no suitable EV is found. In our simulation, the following three selection strategies are employed:
\begin{itemize}
    \item Power-of-$d$ vehicles dispatch
    \item Closest vehicle dispatch
    \item Closest available vehicle dispatch
\end{itemize}

\noindent We explain these algorithms in detail below. Additionally, our simulation framework is modular. It supports testing any possible combination of these first and second decision options.

\subsubsection{Power-of-$d$ Vehicles Dispatch Algorithm}
\noindent After getting a list of available vehicles based on the first decision, the power-of-$d$ vehicles dispatch algorithm optimized the dispatch of EVs by considering the $d \in \mathbb{R}_+$ closest electric vehicles to a customer arrival from the list, where d is a positive integer. Among these $d$ vehicles, the algorithm selects the one with the highest SoC. If all the vehicles are busy or if none of the $d$ closest vehicles have sufficient SoC to serve the trip, the trip is dropped, and the customer leaves the system immediately.

\begin{algorithm}[H]
\caption{Power-of-$d$ Vehicles Dispatch} \label{alg:powerd}
{\fontsize{10}{13}\selectfont
\SetAlgoLined
\SetKwInOut{Input}{Input}
\SetKwInOut{Initialize}{Initialization}
\Input{Trip, List of available EVs, Fix $d \geq 1$}
\For{Every Arrival}{
\If{All vehicles are busy}{
The arrival is dropped\;
}
Choose the $d$ closest vehicles from the list of available EVs\;
\If{None of the $d$ vehicles have enough SoC for the trip}{
The arrival is dropped\;
}
Match the arrival to the vehicle with the highest SoC among the $d$ closest\;
}
}
\end{algorithm}

\noindent The power-of-$d$ vehicles dispatch algorithm can be operated with either a fixed value of $d$ or an adaptive value of $d$. With a fixed $d$, the value of $d$ remains constant throughout the entire simulation period. With an adaptive $d$, the algorithm dynamically adjusts the value of $d$ based on real-time conditions, specifically the number of idle cars with a high SOC. Starting with a predefined value of $d$, the system periodically calculates the average number of such idle cars and increments $d$ if this average exceeds a predefined threshold. Conversely, if the average falls to zero and $d$ is greater than 1, it decrements $d$. This adjustment occurs every specified number of trips to ensure the system remains responsive to changes in vehicle availability and usage patterns, thereby optimizing dispatch efficiency.

\subsubsection{Closest Vehicle Dispatch Algorithm}
\noindent The closest dispatch algorithm is a special case of the power-of-$d$ dispatch policy, where $d=1$. This algorithm first identifies the closest EV to a customer arrival from the available EV list. Once the closest EV is determined, the algorithm checks its SoC. If the vehicle's SoC is above a predefined threshold, which means it has sufficient charge to complete the trip, the vehicle is dispatched to the customer. However, if the SoC is below the threshold, the trip will be dropped, and the customer will leave the system without being served.

\subsubsection{Closest Available Vehicle Dispatch Algorithm}
\noindent With a list of available vehicles from the first decision, the closest available dispatch algorithm also selects the EV from the list. However, it improves on the basis of the simplest closest dispatch algorithm by considering not just the closest EV but also its capability to serve a trip according to its SoC. When a customer arrives, the algorithm considers the closest EV to a customer arrival that has enough SoC to serve the trip. In other words, it initially identifies the closest EV to the pickup location. If the SoC of the current closest EV is too low to complete the trip, the algorithm then proceeds to the next closest EV. This process repeats until it finds one with a sufficient SoC to dispatch. If no suitable EV is found at the end, the trip request is dropped, and the customer leaves the system.

\begin{algorithm}[H]
\caption{Closest Available Dispatch Algorithm} \label{alg:closest_available_dispatch}
{\fontsize{10}{13}\selectfont
\SetAlgoLined
\SetKwInOut{Input}{Input}
\SetKwInOut{Initialize}{Initialization}
\Input{Trip, List of available EVs}
\For{Every Arrival}{
Choose the list of vehicles from the list of available EVs\;
Choose the vehicles that have enough SoC for the trip and sort them based on pickup time\;
\If{None of the vehicles have enough SoC for the trip}{
The arrival is dropped\;
}
Match the arrival to the vehicle with the shortest pickup time and the highest SoC\;
}
}
\end{algorithm}

\noindent Below are examples of the three matching algorithms applied to the same scenario. In these examples, we define the available EV list to include only the cars that are currently idle, charging, or waiting at the charger. Therefore, cars that are driving with passengers or picking up passengers in the diagrams are not included in the available EV list.

\begin{figure}[H]
    \centering
    \includegraphics[scale=0.26]{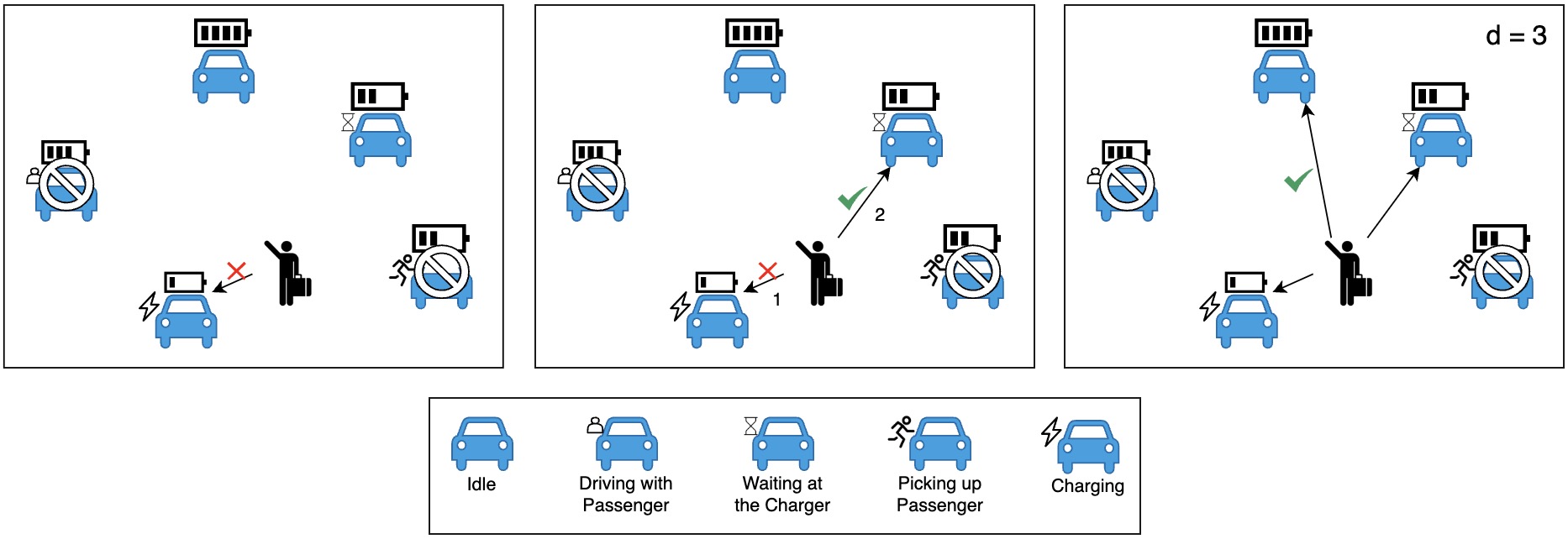}
    \captionsetup{justification=centering, margin=3cm}
    \vspace{-0.3cm}
    \caption{Matching Algorithms Diagrams: (left) Closest Dispatch, (middle) Closest Available Dispatch, (right) Power-of-$d$ ($d$ = 3)}
    \label{fig:matching_diagrams}
\end{figure}

\noindent \textbf{Example 1}: In the leftmost diagram of Fig. \ref{fig:matching_diagrams}, among the available cars, the closest dispatch algorithm chooses the closest EV to the customer arrival. The algorithm determines the car that is currently charging is the nearest car, but it realizes that this car is not able to serve the customer with the current SoC. Therefore, the trip is finally dropped, and the customer is not served.

\noindent \textbf{Example 2}: In the middle diagram of Fig. \ref{fig:matching_diagrams}, from the available EV list, the closest available dispatch algorithm intends to select the closest car to the customer arrival with enough SoC to serve the trip. Initially, it checks the nearest EV to the pickup location, the one that is currently charging, but then it discovers that the car does not have sufficient battery to complete the trip. Consequently, the algorithm moves to the next nearest EV, which is the waiting at the charger one. After determining that the SoC of this car is adequate, it matches this car with the passenger, and the customer is served at the end.

\noindent \textbf{Example 3}: In the rightmost diagram of Fig. \ref{fig:matching_diagrams}, the power-of-$d$ algorithm identifies the $d$ closest cars to the customer's arrival location from the available EV list. Since $d$ is set to $3$ in this case, the cars in the diagram that are idle, waiting at the charger, and charging are being considered. Among these $d$ vehicles, the algorithm chooses the one with the highest SoC, which is the EV that is idle. This car is then matched with the passenger, and the customer is eventually served.

\noindent As a result, although the scenario is exactly the same, the three different matching algorithms will pick three different cars at the end.

\subsection{Charging Algorithms}
\noindent The charging algorithm is prompted every time there is a customer arrival in the system. It is crucial for managing the EV fleets, and it is defined based on a vector of thresholds, $C(t)$, for each time $t$ during the simulation. The charging algorithm has three main components: (1) Deciding which EVs to send to charge; (2) Identifying which chargers are available; (3) Matching the EVs to the chargers and deciding charge duration.

\subsubsection{First Part: Identification of EVs for Charging}

\noindent We define $T$ as the total duration of the simulation. For all $t \in [0, T]$, a threshold $C(t) \in [0, 1]$ is set. At each given time $t$, the algorithm fetches all EVs that are idle. If the EV satisfies $\operatorname{SoC}(t) \leq C(t)$, it is selected for sending to a charger.
\noindent In our simulation, the charging strategy is implemented with two values of thresholds:
\begin{itemize}
    \item Charge idle vehicles all the time: set $C(t) = 0.95$ for all $t \in [0, T]$. This approach promotes frequent charging, implying that the system sends all EVs to charge as soon as they finish serving a trip, as long as their final SoC is less than 95\%. This represents a greedy charging policy.
    \item Charge idle vehicles primarily at night: set $C(t) = 0.4$ during daytime hours and $C(t)=0.95$ in the nighttime. In particular, we define daytime as 6 am - 11 pm, based on the peaks in our customer arrival dataset, and consider the remainder of the time as nighttime. This strategy minimizes daytime charging with a low threshold, and the vehicles are only sent to charge if absolutely necessary, allowing them to be able to serve more trips when demand is higher. By setting a high threshold at night, the algorithm also ensures that the vehicles are charged during off-peak hours when they are less likely to be in use.
\end{itemize}

\subsubsection{Second Part: Identification of Availability of Chargers}

\noindent The system first selects a subset of chargers and determines their availability based on the number of available posts and the number of EVs driving to each charger. Specifically, a charger is deemed available if the number of available posts exceeds the product of $\alpha$ and the number of EVs driving to that charger. Here, $\alpha$ is a parameter that can be set by the user.

\noindent We observe a trade-off between maximizing utilization of the chargers and minimizing the delay experienced by the EVs at the chargers depending on the value of $\alpha \in [0, 1]$. If we set $\alpha=1$, we assume that the driving EV is already at the charger while it is actually still on its way to the charger, which means we are wasting the capacity of the charger. If we set $\alpha=0$, we might end up sending too many EVs to the same charger. Since there is a delay in them reaching the charger, they could all arrive simultaneously, leading to excessive numbers of EVs waiting at the charger. In our simulation, we test the system's performance for $\alpha \in \{0, 0.5, 1\}$ to optimize charger allocation.

\subsubsection{Third Part: Matching EVs to Chargers}

\noindent Upon determining the available chargers and the EVs that need to be sent to the chargers, we proceed to match each EV with an appropriate charger. The matching process involves two policies:
\begin{itemize}
    \item Closest Available Charger: This approach prioritizes proximity. It matches each EV to the nearest available charger to minimize the travel time from the EV's current location to the charger location. The algorithm of Closest Available Charger is similar to Algorithm \ref{alg:closest_available_dispatch}.
    \item Power-of-$d$ Choices for Charging: This way evaluates the $d$ closest chargers to each EV. If all these chargers are unavailable, as defined in the second step, the EV is not sent to charge. The algorithm of Power-of-$d$ Choices for Charging is similar to Algorithm \ref{alg:powerd}.
\end{itemize}
Once an EV is sent to charge, it is charged up to 100\% unless the process is interrupted by Algorithm \ref{alg: interrupt_charging}.


\section{Running the Simulation}

\subsection{Real Life Data Input}
\noindent To create a realistic simulation framework, we aim to capture various real-life phenomena, such as peaks and valleys in the trip demand, congestion on the streets, and inhomogeneous charging patterns (e.g., higher charging demand during nighttime). The sample simulation results utilize open-source NYC taxi data, which is the May 2024 TLC Yellow Taxi Trip Record Data \cite{nyc2023}. This dataset records trips taken by passengers using the Yellow Cab in New York City. Thus, capturing the inhomogeneous demand patterns and varying trip velocities due to congestion. Our simulation framework is also compatible with various open-source datasets available out there, such as RideAustin \cite{rideaustin2023} and Chicago taxi data \cite{chicago2023}.

\noindent We use the NYC taxi dataset to extract key trip information, such as the origin, destination, arrival time, and trip distance of all trips. During the data cleaning process, we filter the data based on a specified time interval and remove any erroneous data entries. For instance, we set the lower and upper percentiles to filter out the invalid latitude and longitude values in the dataset.

\begin{wrapfigure}[15]{r}{0.31\textwidth}
    \centering
    \includegraphics[width=0.31\textwidth]{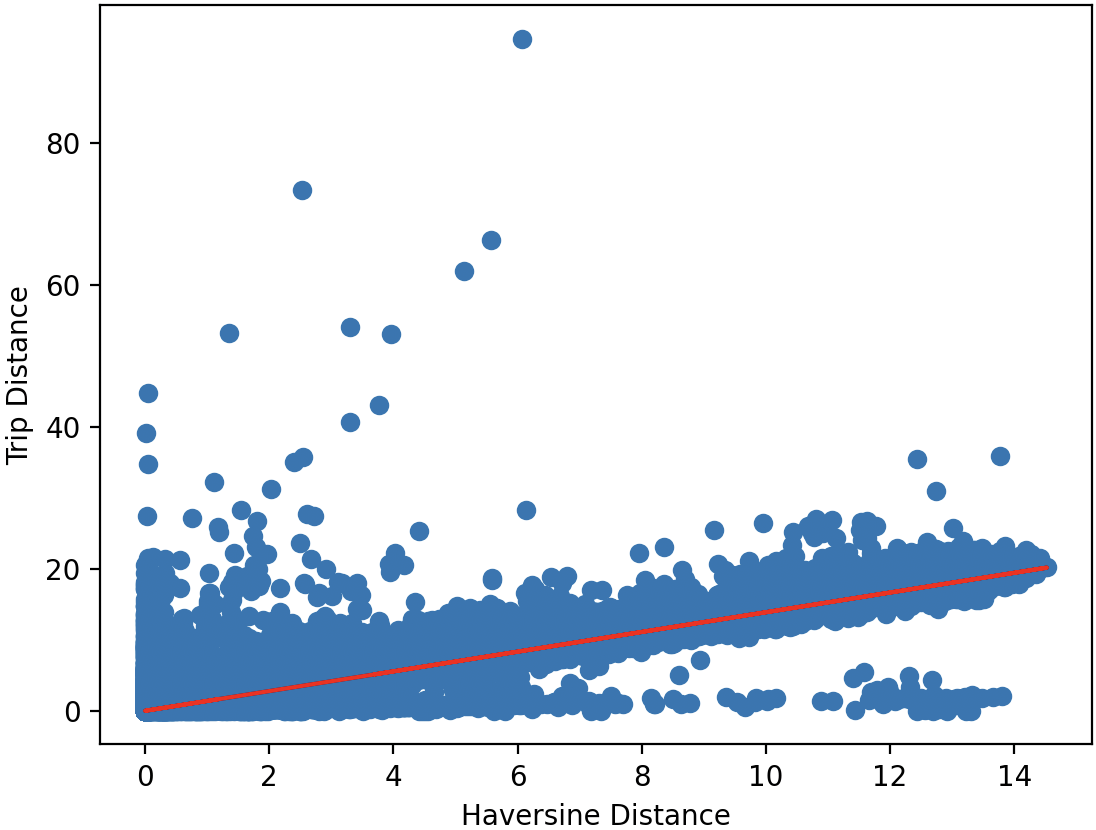}
    \vspace{-0.6cm}
    \caption{Linear Regression}
    \label{fig:linear_regression_main}
\end{wrapfigure}

\noindent Since the pickup time and distance of the trips are not specified in the dataset, we calculate the Haversine distance between pickup and drop-off locations and adjust this using linear regression to better capture the actual geometry of the city. To be specific, we split the sample data into two sets: 80\% of the data for training the model and 20\% for testing the accuracy of the model using mean squared error. After generating the scatter plot, as shown in Fig.~\ref{fig:linear_regression_main}, we fit a linear regression model to the data and treat the slope of the linear regression function as the correction factor. By multiplying the Haversine distance by this correction factor, we obtain a more accurate estimate of the travel distance between any two points in the range of NYC. Based on these distances, we estimate the pickup time for each trip. The R-squared value of the linear regression model is 0.848, with a Mean Squared Error of 1.22, indicating a good fit to the sample trip distances. The framework also supports Manhattan distance calculations by applying the correction factor to the Manhattan distance for greater accuracy.

\subsection{Inputs to the Simulation}

\textbf{Fixed Parameters}: The following parameters remain constant across all sample simulations.

\begin{table}[h]
    \centering
    \begin{tabular}{ |p{5.5cm}|p{8cm}| }
        \hline
        Dataset & New York Taxi: May 1, 2024 to May 4, 2024\\
        \hline
        EV Battery Pack Size & 51.25 kWh (after a 10\% battery degredation)\\  
        \hline
        EV Consumption Rate & 230 Wh/mi\\  
        \hline
        EV Charge Rate & 20 kW\\
        \hline
    \end{tabular}
    \caption{Fixed Parameters}
    \label{table:fixed_parameters}
\end{table}

\noindent The input dataset includes trips from the NYC dataset between May 1, 2024, and May 4, 2024. The EV parameters are based on the Tesla Model 3 standard range, with a battery pack size of 51.25 kWh (accounting for a nominal 10\% degradation) and an energy consumption rate of 230 Wh/mi. We assume a charge rate of 20 kW, which is typical for a Type-2 charger. The average velocity of the EVs is calculated from the NYC trip data, which averages approximately 11.21 mi/hr.

\noindent The number of EVs in the simulation is set to match peak daytime trip demand to ensure high throughput in the simulations. The number of chargers is also set proportionally to ensure the chargers can support the charging requirements of the EVs. The initial locations of EVs are chosen randomly from the pickup locations of the trips in the dataset. The locations of charging stations are chosen randomly and uniformly across the NYC range. Each charging station has 4 charging posts.

\section{Simulation Results}

\subsection{Sample Simulation Outcomes}

The fleet size is set to match the peak demand, which is 2101, according to the simulation data. The simulation runs over three days and applies the power-of-$d$ vehicle dispatch algorithm with $d = 10$. The charging strategy continuously charges idle vehicles throughout the day. The average speed of EVs is 11.21 mph, which is derived from the dataset. The number of charging stations is calculated as the proportion of the fleet needing to charge to compensate for the consumption due to driving, resulting in 270 chargers.

\begin{figure}[tbh!]
    \hspace{-2cm}
    \centering
    \begin{minipage}[h]{0.45\textwidth}
        \centering
        \includegraphics[scale=0.25]{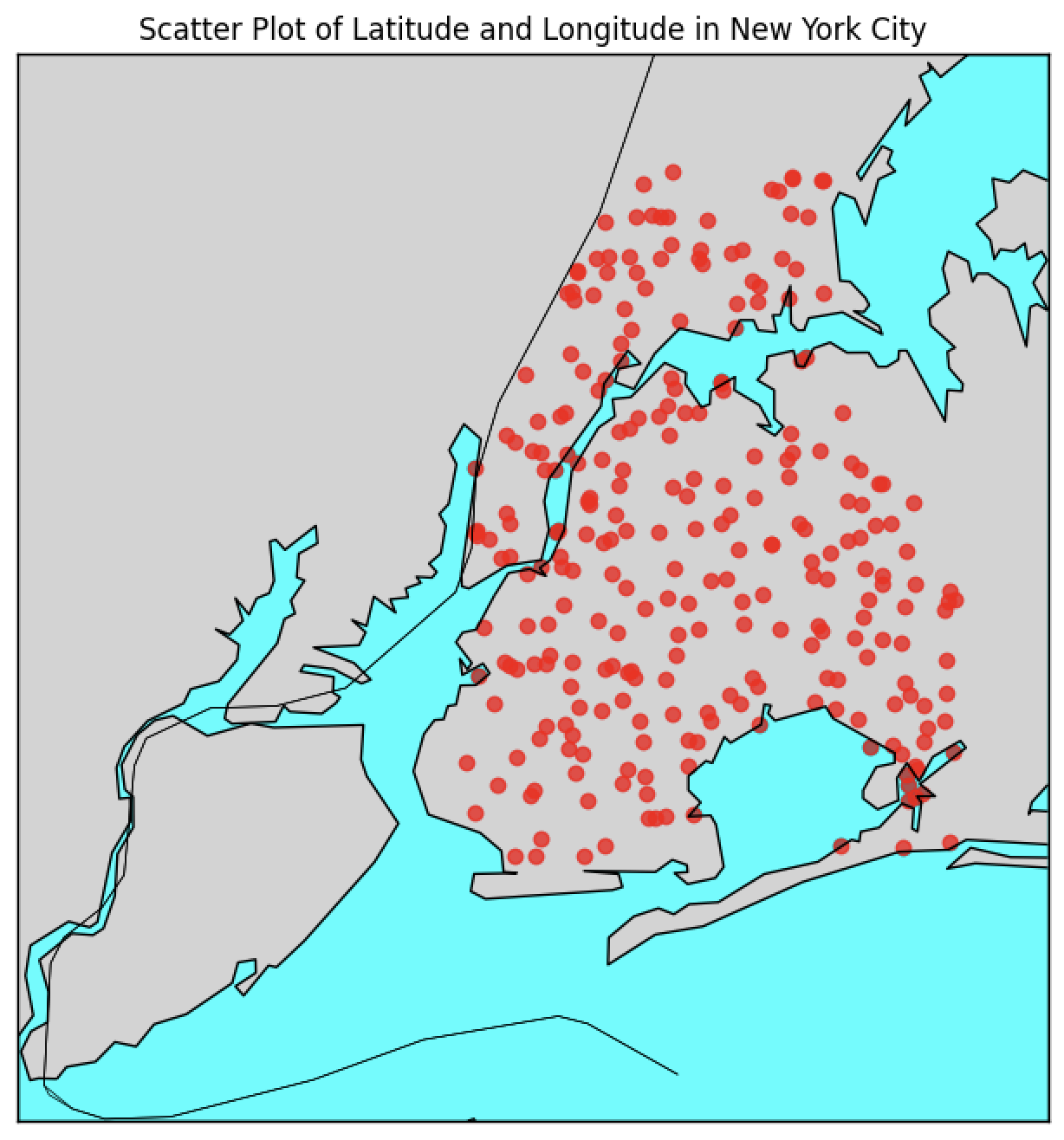}
        \captionsetup{justification=centering, margin=1cm}
        \caption{Charger locations of current simulation}
        \label{fig:charger_locations}
    \end{minipage}
    \hspace{-0.5cm}
    \begin{minipage}[h]{0.45\textwidth}
        \centering
        \includegraphics[scale=0.5]{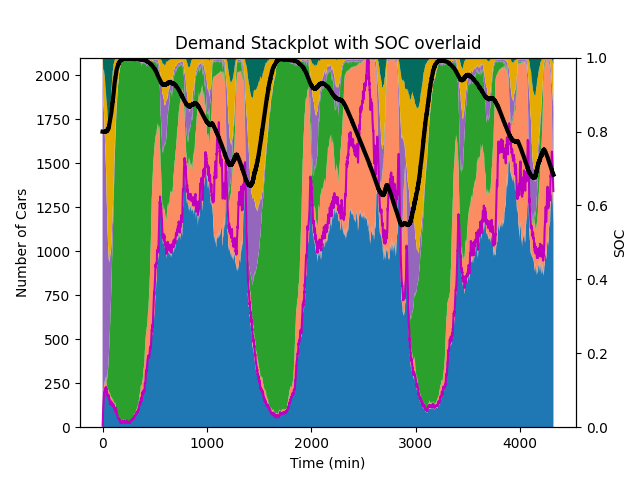}
        \includegraphics[scale=0.21]{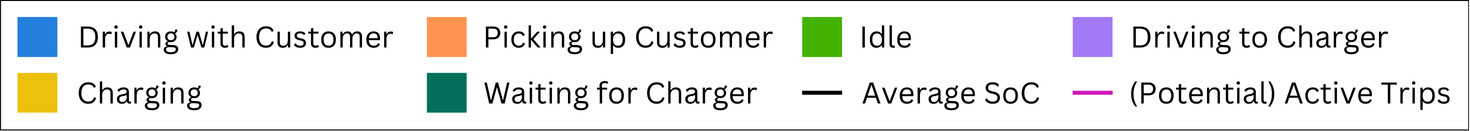}
        \captionsetup{justification=centering, margin=1cm}
        \vspace{-0.5cm}
        \caption{3-day power-of-$10$ and charge idle cars all the time}
        \label{fig:3_day_d_10_low_velocity_n_chargers}
    \end{minipage}
\end{figure}

\begin{figure}[tbh!]
    \hspace{-2cm}
    \centering
    \begin{minipage}[h]{0.45\textwidth}
        \centering
        \includegraphics[scale=0.45]{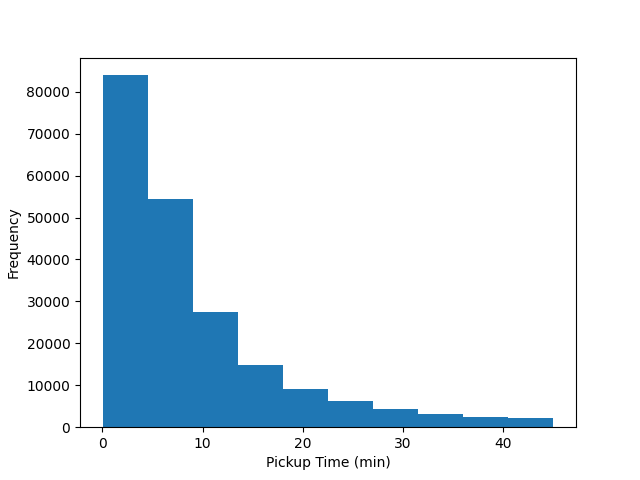}
        \captionsetup{justification=centering, margin=0.1cm}
        \vspace{-0.2cm}
        \caption{Pickup Time Frequency}
        \label{fig:pickup_time_hist}
    \end{minipage}
    \hspace{-0.5cm}
    \begin{minipage}[h]{0.45\textwidth}
        \centering
        \begin{table}[H]
            \begin{tabular}{ |p{4.2cm}|p{4.2cm}| }
                \hline
                \multicolumn{2}{|c|}{3-day power-of-$10$ and charge idle cars all the time} \\
                \hline
                \hline
                Matching algorithm & power-of-$10$\\
                \hline
                Charging algorithm & charge idle cars all the time\\
                \hline
                Average velocity (mph) & 11.21\\
                \hline
                Fleet size & 2101\\
                \hline
                Number of chargers & 270\\
                \hline
                Average trip time (mins) & 17.29\\
                \hline
                Average trip distance (miles) & 3.32\\
                \hline
                Average pickup time (mins) & 8.90\\
                \hline
                Average time to chargers (mins) & 10.32\\
                \hline
                Number of trips to charger per car per hour & 0.25\\
                \hline
                Average SoC over time over cars & 0.68\\
                \hline
                Service level & 92.21\%\\
                \hline
                Workload served & 87.05\%\\
                \hline
            \end{tabular}
            \caption{Data Table}
            \label{table:data_table}
        \end{table}
    \end{minipage}
\end{figure}


\noindent Running the simulation generates Fig.~\ref{fig:charger_locations}, \ref{fig:3_day_d_10_low_velocity_n_chargers}, \ref{fig:pickup_time_hist}, and Table \ref{table:data_table}, which collectively illustrate the performance of the implemented matching and charging algorithms. We now explain these results in detail below. 

\vspace{0.2cm}
\noindent \textbf{Charger locations scatter plot (Fig.~\ref{fig:charger_locations})}: The location of each charger is randomly selected from a sample of pickup locations. During the simulation, we filter out pickup locations by refining the percentile ranges for latitude and longitude to ensure that chargers are placed within specified bounds, avoiding locations with minimal trip demand. The scatter plot displays the distribution of charger locations across New York City, with each red dot representing an individual charger. The chargers are spread out across the city, with a concentration in Manhattan, Brooklyn, and parts of Queens. Since the charger distribution is based on pickup locations, it aligns with areas of higher demand for electric vehicles and greater population density.

\vspace{0.2cm}
\noindent \textbf{Demand curve stack plot (Fig.~\ref{fig:3_day_d_10_low_velocity_n_chargers})}: The left vertical axis represents the number of vehicles, while the right vertical axis measures the average SoC of all EVs, and the horizontal axis shows time in minutes over a three-day period. The plot breaks down EV distribution across various activities. The majority of vehicles are either idle (light green area) or serving passengers (blue area), with smaller portions of the fleet traveling to pickup locations (orange area), charging (yellow area), driving to charging stations (purple area), or waiting for a charger (dark green area). The black curve, representing the average SoC of all EVs over time, follows a cyclical pattern, decreasing when vehicles are in use and increasing as they charge, which indicates typical daily usage and charging cycles. When the average SoC drops, the yellow area corresponding to charging activity begins to rise, and conversely, it decreases when SoC levels recover. The large light green area suggests a significant number of idle vehicles during nighttime. The pink line indicates potential active trips over time, and the extensive blue area beneath this line shows that a high percentage of these trips are successfully served. The small purple areas suggest that vehicles spend relatively little time driving to charging stations, indicating efficiency in such process. Overall, the power-of-$10$ dispatch algorithm appears to effectively manage the fleet, ensuring high passenger service rates while also maintaining sufficient SoC levels through consistent charging of idle vehicles throughout the day. 

\vspace{0.2cm}
\noindent \textbf{Pickup time frequency distribution (Fig.~\ref{fig:pickup_time_hist})}: The histogram presents the frequency of electric vehicle pickup times over the three-day simulation. Most pickups occur within the first 10 minutes. Thereafter, the frequency of pickups drops significantly, indicating rapid and efficient service at shorter time intervals. The sharp decline in pickup times highlights the system's effectiveness in achieving quick pickups and minimizing passenger wait times.

\vspace{0.2cm}
\noindent \textbf{Data table (Tab.~\ref{table:data_table})}: The table provides key performance metrics. The average trip time is 17.29 minutes, and the average trip distance is 3.32 miles. The average pickup time indicates the mean time for EVs to reach passengers, and the average drive time to the charger is the average duration for EVs to travel to the charging stations. Both the average pickup time of 8.90 minutes and the average time to reach chargers of 10.32 minutes are relatively low, suggesting that the matching algorithm and charging strategy are effectively minimizing passenger wait times and maintaining operational efficiency. The number of trips to a charger per car per hour refers to how often each vehicle needs to visit a charging station within an hour. With a rate of 0.25, each EV, on average, requires charging roughly once every four hours. This low frequency suggests that the vehicles have sufficient range and that the charging strategy effectively minimizes service disruptions. The average SoC over time across all cars is 0.68, reflecting a relatively high SoC that ensures vehicles maintain a healthy charge, allowing them to remain readily available for passenger service. The service level denotes the percentage of total trips that are completed, and the workload refers to the percentage of total miles that are served. The service level of 92.21\% and workload of 87.05\% are both relatively high, which reflect the simulation's overall strong performance.

\subsection{Comparative Analysis of Simulations with Different Parameters:}

\subsubsection{Scenario 1: Lower Charger Parameter}

\noindent This scenario uses the fixed parameters outlined in Table \ref{table:fixed_parameters}. 
The fleet size is the same as peak trip demand. The average velocity of the EVs is 11.21 mph, which is derived from the dataset. 
Different simulations will employ various matching and charging algorithms as specified.

\begin{table}[h]
    \centering
    \vspace{-0.2cm}
    \begin{tabular}{ |p{2cm}|p{1.65cm}|p{1.65cm}|p{1.65cm}|p{1.65cm}|p{1.65cm}|p{1.65cm}|p{1.65cm}| }
        \hline
        Matching  & closest available dispatch & power-of-$1$ (closest dispatch) & power-of-$5$ & power-of-$10$ & power-of-$d$ adaptive & power-of-$10$ & power-of-$10$\\
        \hline
        Charging  & \multicolumn{5}{c|}{charge idle cars all the time}  & \multicolumn{2}{c|}{charge primarily at night} \\
        \hline
        Interrupt Charging & yes & yes & yes & yes & yes & yes & no\\
        \hline
        Average pickup time (mins) & 8.14 & 5.20 & 7.20 & 8.12 & 8.10 & 7.18 & 9.81\\
        \hline
        Average time to chargers (mins) & 23.67 & 21.65 & 28.07 & 27.35 & 26.80 & 18.69 & 30.19\\
        \hline
        Service level & 88.48\% & 79.67\% & 89.63\% & 91.12\% & 91.54\% & 91.58\% & 88.52\%\\
        \hline
        Workload served & 83.55\% & 76.52\% & 85.34\% & 86.27\% & 86.41\% & 86.51\% & 83.73\%\\
        \hline
    \end{tabular}
    \caption{3-day Simulations with average velocity of 11.21 mi/hr, fleet size of 2101, and a total of 67 chargers.}
    \label{table:3_day_low_velocity_low_charger}
\end{table}

    

\begin{figure}[bth!]
    \centering
    \includegraphics[scale=0.25]{figures/Legend.png}
    \begin{minipage}[bth!]{0.45\textwidth}
        \centering
        \includegraphics[scale=0.45]{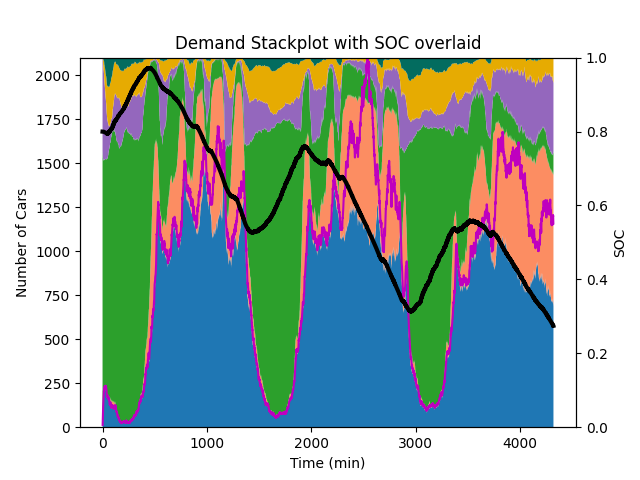}
        \captionsetup{justification=centering, margin=0.5cm}
        \vspace{-0.3cm}
        \caption{3-day closest available dispatch and charge all the time with v=11.21mph and less chargers}
        \label{fig:3_day_closest_available_low_velocity_new}
    \end{minipage}
    \hspace{1cm}
    \begin{minipage}[bth!]{0.45\textwidth}
        \centering
        \includegraphics[scale=0.45]{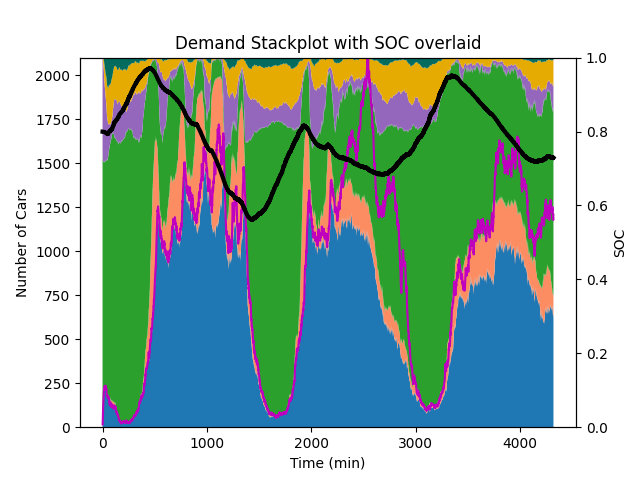}
        \captionsetup{justification=centering, margin=0.5cm}
        \vspace{-0.3cm}
        \caption{3-day power-of-$1$ (closest dispatch) and charge all the time with v=11.21mph and less chargers}
        \label{fig:3_day_d_1_low_velocity_n_chargers_0.25_new}
    \end{minipage}
\end{figure}
\begin{figure}[bth!]
    \centering
    \begin{minipage}[h]{0.32\textwidth}
        \centering
        \includegraphics[scale=0.33]{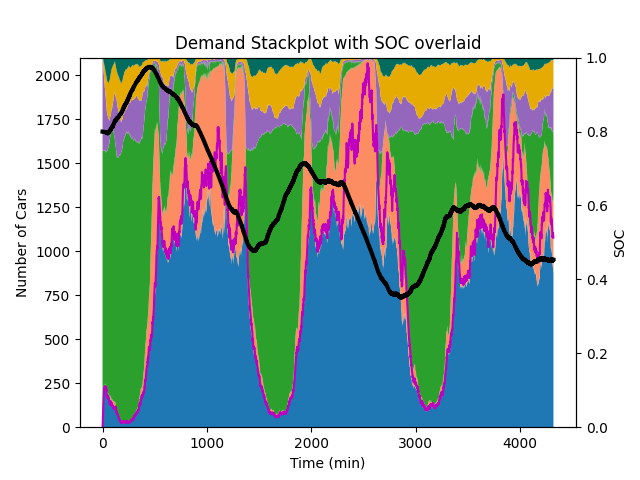}
        \captionsetup{justification=centering, margin=0.6cm}
        \vspace{-0.3cm}
        \caption{3-day power-of-$10$ and charge all the time with v=11.21mph and less chargers}
        \label{fig:3_day_d_10_low_velocity_n_chargers_0.25_new}
    \end{minipage}
    \hfill
    \begin{minipage}[bth!]{0.32\textwidth}
        \centering
        \includegraphics[scale=0.33]{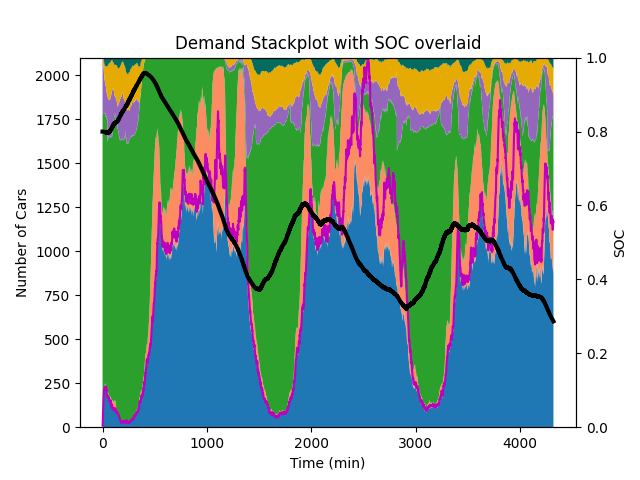}
        \captionsetup{justification=centering, margin=0.1cm}
        \vspace{-0.3cm}
        \caption{3-day power-of-$10$, charge primarily at night with v=11.21mph and less chargers, interrupt charging}
        \label{fig:3_day_d_10_low_velocity_n_chargers_0.25_night_new}
    \end{minipage}
    \hfill
    \begin{minipage}[bth!]{0.32\textwidth}
        \centering
        \includegraphics[scale=0.33]{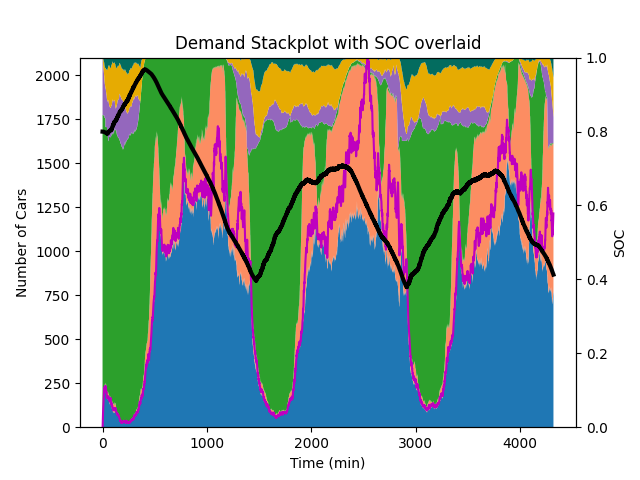}
        \captionsetup{justification=centering, margin=0.1cm}
        \vspace{-0.3cm}
        \caption{3-day power-of-$10$ and charge primarily at night with v=11.21mph and less chargers, no interrupt charging}
        \label{fig:3_day_d_10_low_velocity_n_chargers_0.25_night_no_interrupt_charging}
    \end{minipage}
\end{figure}

\noindent \textbf{Closest Available Vehicle Dispatch vs. Power-of-$d$ Vehicles Dispatch with Fixed $d=1$:}

\noindent In comparing the simulation using the closest available vehicle dispatch with the simulation employing power-of-$1$ (closest vehicle dispatch), Fig. \ref{fig:3_day_closest_available_low_velocity_new} exhibits a prominent orange region during the peak demand hours, indicating a significant number of vehicles engaged in pickups, whereas Fig. \ref{fig:3_day_d_1_low_velocity_n_chargers_0.25_new} displays a large green area representing idle vehicles instead in the power-of-$1$ algorithm. This difference arises because the closest available dispatch algorithm selects the nearest vehicle that has sufficient SoC to complete the trip, and it is continuously searching until such a vehicle is found, leading to a higher average pickup time and more vehicles actively engaged in pickups (larger orange area) in Fig. \ref{fig:3_day_closest_available_low_velocity_new}. On the other hand, the closest dispatch algorithm immediately drops trip requests if the nearest vehicle lacks sufficient SoC, thereby increasing vehicle idle time (green area). Despite the longer pickup times, the closest available dispatch algorithm achieves a higher service level and completes more trips than the closest dispatch algorithm. From another perspective, since the closest dispatch algorithm, which frequently drops trips, results in many vehicles remaining idle during peak hours, its average SoC recovers quickly. However, the SoC curve for the closest available dispatch algorithm shows a substantial decline.

\vspace{0.2cm}
\noindent \textbf{Closest Available Vehicle Dispatch vs. Power-of-$d$ Vehicles Dispatch with Fixed $d=10$:}

\noindent When comparing the closest available vehicle dispatch algorithm to the power-of-$10$ vehicle dispatch algorithm, Fig. \ref{fig:3_day_closest_available_low_velocity_new} shows a significantly larger orange area on the third day than Fig. \ref{fig:3_day_d_10_low_velocity_n_chargers_0.25_new}. The closest available dispatch algorithm, due to its longer pickup times, results in a significant drop in the average SoC, which limits its capacity to serve many trips on day 3. In contrast, the power-of-$10$ dispatch algorithm stabilizes the SoC curve so that it is able to handle many more trips on the third day. This is evident in Fig. \ref{fig:3_day_d_10_low_velocity_n_chargers_0.25_new}, where the blue area (served trips) is larger, reflecting better overall performance. The power-of-$10$ algorithm achieves higher service levels and workload efficiency by balancing vehicle utilization and maintaining a more consistent SoC across the fleet.



\vspace{0.2cm}
\noindent \textbf{Power-of-$d$ Vehicle Dispatch with Fixed $d=10$ vs. Power-of-$d$ Vehicles Dispatch with Adaptive $d$:}

\noindent The adaptive power-of-$d$ dispatch algorithm begins with an initial $d$ value of $5$ and dynamically adjusts this value throughout the simulation based on the average number of idle vehicles with high SoC. Starting with $d=5$, the system periodically calculates the average number of idle vehicles and increases $d$ if this average exceeds a predefined threshold (set to $5\%$) while incoming trips are dropped. Conversely, if the average number of idle vehicles drops to zero and $d$ is greater than $1$, the algorithm decrements $d$. This evaluation is performed every 1000 trips, ensuring the system remains responsive to changes in vehicle availability and usage patterns, thereby optimizing dispatch efficiency. These adjustments are aimed at optimizing dispatch efficiency by considering the most suitable number of vehicles under varying system conditions and performance metrics. During the simulation, $d$ is gradually increased or decreased as needed, and the peak value that it has reached is 15. The adaptive power-of-$d$ approach results in an average pickup time that lies between those of the fixed power-of-$5$ and power-of-$10$ algorithms. The balance between flexibility and efficiency enables the adaptive power-of-$d$ algorithm to outperform not only the closest available dispatch algorithm but also the fixed power-of-$d$ algorithms with $d=1, 5, 10$. Consequently, the adaptive $d$ algorithm consistently delivers the highest service levels and workload efficiency among all matching algorithms tested in Table \ref{table:3_day_low_velocity_low_charger}. Such an adaptive algorithm goes beyond the analysis in \cite{varma_ev} and is our contribution. 

\vspace{0.2cm}
\noindent \textbf{Charge Idle Vehicles All the Time vs. Charge Idle Vehicles Primarily at Night:}

\noindent To compare the ``charge all the time'' algorithm with the ``charge primarily at night'' algorithm, we examine Fig. \ref{fig:3_day_d_10_low_velocity_n_chargers_0.25_new} and Fig. \ref{fig:3_day_d_10_low_velocity_n_chargers_0.25_night_new}, both using the power-of-$10$ vehicle dispatch algorithm. The ``charge all the time'' algorithm directs vehicles to charging stations throughout the day and night, resulting in a continuous yellow area representing charging vehicles in Fig. \ref{fig:3_day_d_10_low_velocity_n_chargers_0.25_new}. This algorithm frequently sends vehicles to chargers, even during peak demand hours, whenever their SoC drops below 0.95. In contrast, the ``charge primarily at night'' algorithm focuses on charging vehicles during the nighttime, between 11 PM and 6 AM. The threshold $C(t)$ is set to $0.4$ during daytime and increases to $0.95$ at night. Vehicles with a SoC below $C(t)$ are directed to charging stations. During the day, vehicles are only sent to chargers if their SoC falls below 0.4, while at night, vehicles are charged more aggressively if their SoC is below 0.95, resulting in the disappearance of the yellow region on the first day in Fig. \ref{fig:3_day_d_10_low_velocity_n_chargers_0.25_night_new}. At 6 AM, the average SoC reaches the local maxima and then begins to decline throughout the daytime. By setting a lower threshold during the daytime, this approach reduces daytime charging, preserving vehicle availability for higher demand periods. By adopting a higher threshold during nighttime, this strategy ensures that charging activities are concentrated during off-peak hours when most vehicles are idle overnight. The "charge primarily at night" algorithm significantly decreases average pickup time and time to chargers, as indicated by the smaller orange and purple areas in Fig. \ref{fig:3_day_d_10_low_velocity_n_chargers_0.25_night_new}. During peak demand on day 2, the yellow region is smaller, and the blue area has increased, suggesting better workload management. Consequently, this nighttime charging strategy enhances overall service level and workload efficiency, making it more effective for meeting high trip demand during daytime hours.

\noindent Furthermore, all simulations, except for Fig. \ref{fig:3_day_d_10_low_velocity_n_chargers_0.25_night_no_interrupt_charging}, were run with an interrupt charging strategy. Comparing the last two columns in Table \ref{table:3_day_low_velocity_low_charger}, represented by Fig. \ref{fig:3_day_d_10_low_velocity_n_chargers_0.25_night_new} and Fig. \ref{fig:3_day_d_10_low_velocity_n_chargers_0.25_night_no_interrupt_charging}, the advantages of interrupt charging become evident. The ``charge primarily at night'' algorithm largely improved the performance when interrupt charging was enabled. This is because the system considers not only the idle vehicles, but also those driving to chargers, charging, and waiting for chargers. By pausing the charging activities for these vehicles when demand peaks, the system optimizes resource allocation. For example, in Fig. \ref{fig:3_day_d_10_low_velocity_n_chargers_0.25_night_new}, far more trip are completed compared to Fig. \ref{fig:3_day_d_10_low_velocity_n_chargers_0.25_night_no_interrupt_charging}, particularly during the peak demand period around time = 2600 minutes. Therefore, this combined approach of nighttime charging with interrupt charging maximizes vehicle availability during peak hours, allowing for more efficient utilization of the fleet and significantly improving performance metrics.

\subsubsection{Scenario 2: Higher Charger Parameter}

\noindent This scenario follows the fixed parameters listed in Table \ref{table:fixed_parameters}. 
The fleet size matches the peak trip demand. The average vehicle velocity is 11.21 mph, as derived from the dataset. 
The charging algorithm employed in this scenario is the ``charge idle cars all the time'' algorithm. Various simulations use different matching algorithms, as specified.

\begin{table}[h]
    \centering
    \vspace{-0.2cm}
    \begin{tabular}{ |p{3.1cm}|p{2.5cm}|p{2.4cm}|p{2.1cm}|p{2.1cm}|p{2.1cm}|}
        \hline
        Matching algorithm & closest available & closest dispatch & power-of-$1.1$ & power-of-$3$ & power-of-$5$\\
        \hline
        Avg pickup time (min) & 6.83 & 6.73 & 6.72 & 7.77 & 8.29\\
        \hline
        Avg time to charger (min) & 10.84 & 10.64 & 10.57 & 10.93 & 11.16\\
        \hline
        Service level & 94.18\% & 93.53\% & 94.01\% & 93.01\% & 92.59\%\\
        \hline
        Workload served & 89.15\% & 88.68\% & 88.96\% & 88.03\% & 87.75\%\\
        \hline
    \end{tabular}
    \caption{3-day Simulations with average velocity 11.21 mi/hr, fleet size of 2101, a total of 270 chargers. The charging algorithm is ``charge idle cars at all times'' as in Section \ref{sec: match_charge}.}
    \label{table:3_day_low_velocity_high_charger}
\end{table}

\begin{figure}[h]
    \centering
    \includegraphics[scale=0.25]{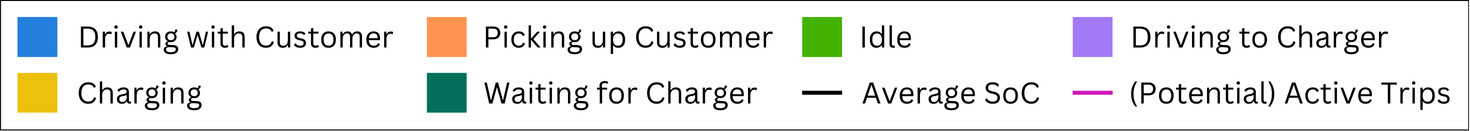}
\end{figure}
\vspace{-0.5cm}
\begin{figure}[H]
    \centering
    \begin{minipage}[h]{0.32\textwidth}
        \centering
        \includegraphics[scale=0.33]{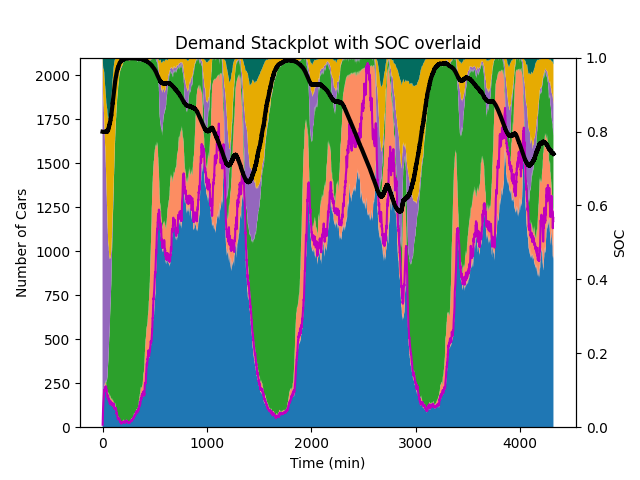}
        \captionsetup{justification=centering, margin=0.1cm}
        \vspace{-0.3cm}
        \caption{3-day power-of-$1$ (closest dispatch) with v=11.21mph and more chargers}
        \label{fig:3_day_d_1_low_velocity_n_chargers_new}
    \end{minipage}
    \hfill
    \begin{minipage}[h]{0.32\textwidth}
        \centering
        \includegraphics[scale=0.33]{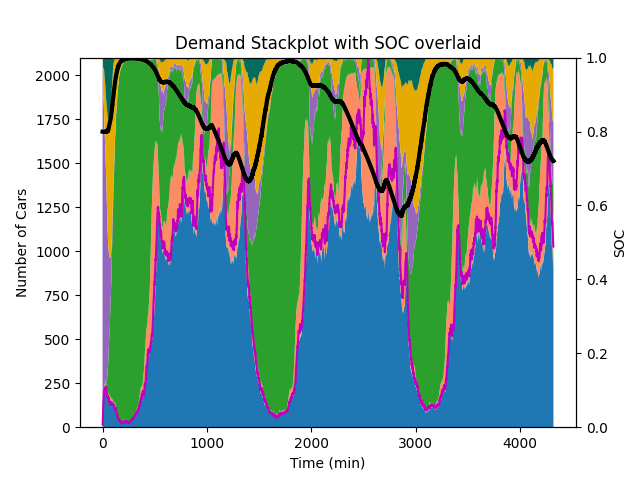}
        \captionsetup{justification=centering, margin=0.3cm}
        \vspace{-0.3cm}
        \caption{3-day power-of-$1.1$ with v=11.21mph and more chargers}
        \label{fig:3_day_d_1.1_low_velocity_n_chargers_new}
    \end{minipage}
    \hfill
    \begin{minipage}[h]{0.32\textwidth}
        \centering
        \includegraphics[scale=0.33]{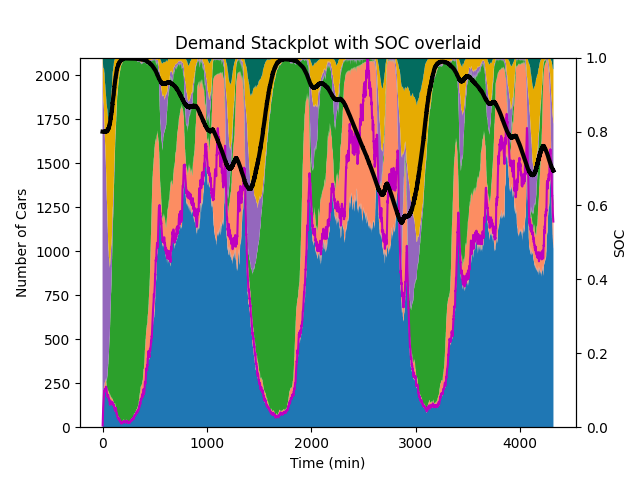}
        \captionsetup{justification=centering, margin=0.5cm}
        \vspace{-0.3cm}
        \caption{3-day power-of-$5$ with v=11.21mph and more chargers}
        \label{fig:3_day_d_5_low_velocity_n_chargers_new}
    \end{minipage}
\end{figure}

\noindent \textbf{Comparison of Power-of-$d$ Vehicle Dispatch Algorithms with Different $d$ in 3-day Simulation:}

\noindent According to Table \ref{table:3_day_low_velocity_high_charger}, the power-of-$1.1$ vehicle dispatch algorithm delivers the highest service level and workload efficiency among all tested power-of-$d$ algorithms, with the shortest average pickup time and the lowest average time spent driving to chargers. When $d$ is not an integer, we consider $\lfloor d \rfloor$ closest vehicles with probability $\lceil d \rceil - d$ and $\lceil d \rceil$ closest vehicles otherwise. And so, we are considering $d$ closest vehicles averaged over time. Such an approach allows us to fine-tune the value of $d$ and it helps improve the system performance.


\noindent Comparing $d=1.1$ (Fig. \ref{fig:3_day_d_1.1_low_velocity_n_chargers_new}) with $d=1$ (Fig. \ref{fig:3_day_d_1_low_velocity_n_chargers_new}) and $d=5$ (Fig. \ref{fig:3_day_d_5_low_velocity_n_chargers_new}), the power-of-$1.1$ outperforms the other two algorithms during peak demand around time = 2500 minutes. The closest dispatch algorithm (power-of-$1$) restricts the consideration to the nearest EV for each customer, and the trips are dropped if the closest vehicle lacks sufficient SoC. Thus, while other nearby EVs might have adequate SoC, we still do not serve the customer. This suggests that increasing the value of $d$ could enhance performance by broadening the selection pool of candidate vehicles. The power-of-$3$ and power-of-$5$ algorithms, with higher values of $d$, evaluate broader selections of vehicles but result in longer pickup times, as indicated by the extensive orange area in Fig. \ref{fig:3_day_d_5_low_velocity_n_chargers_new}, which eventually results in a smaller proportion of customers served.

\noindent These findings demonstrate that, under conditions of higher charger availability, the power-of-$1.1$ vehicle dispatch algorithm offers the highest service level and workload efficiency among the four tested power-of-$d$ algorithms. The value of $d=1.1$ is the optimal value of $d$ among the four simulations since it effectively minimizes pickup times while also balancing the load across the fleet.

\noindent Determining the optimal value of $d$ involves evaluating the trade-off between minimizing pickup times and balancing the load across the fleet. In Fig. \ref{fig:3_day_d_5_low_velocity_n_chargers_new}, a lower $d$ can effectively decrease the average pickup time, enhancing overall performance. In Fig. \ref{fig:3_day_d_1_low_velocity_n_chargers_0.25_new}, increasing $d$ can reduce the excessive number of idle vehicles during peak demand periods. In Fig. \ref{fig:3_day_d_1.1_low_velocity_n_chargers_new}, an optimal value of $d$ yields the best service level and workload efficiency, enabling the system to serve a higher number of trips during peak demand times. Therefore, an optimal dispatch strategy would aim to minimize idle vehicles while efficiently handling trip demands without spending too much time picking up customers.

\vspace{0.2cm}
\noindent \textbf{Lower Charger Parameter vs. Higher Charger Parameter:}

\noindent Comparing Table \ref{table:3_day_low_velocity_low_charger} (low charger availability) and Table \ref{table:3_day_low_velocity_high_charger} (high charger availability), the latter shows a significant decrease in the average travel time for vehicles to reach chargers. Additionally, when all other parameters are held constant, comparing Fig. \ref{fig:3_day_d_1_low_velocity_n_chargers_0.25_new} with fewer chargers to Fig. \ref{fig:3_day_d_1_low_velocity_n_chargers_new} with more chargers, the latter scenario, with increased charger availability, substantially raises the average SoC. The number of vehicles driving to chargers at any given time largely decreases, as indicated by the much smaller purple area in Fig. \ref{fig:3_day_d_1_low_velocity_n_chargers_new} compared to Fig. \ref{fig:3_day_d_1_low_velocity_n_chargers_0.25_new}. Hence, it becomes easier for the SoC curve to revert to its peak value when there are more chargers available. This ensures that vehicles start each operational cycle with higher energy reserves, reducing the risk of trip drops due to low SoC.

\section{Conclusions and Future Work}

\subsection{Conclusions}

\noindent This paper presented a simulation-based analysis of electric vehicle fleet management, examining the performance of various matching and charging algorithms. By evaluating these different dispatch strategies and charging algorithms across multiple scenarios, the study provided insights into their effects on key performance metrics, including service levels, workload efficiency, average time picking up passengers, average time driving to chargers, and overall fleet efficiency.

\noindent The simulations revealed that power-of-$d$ dispatch strategies are more efficient in improving fleet utilization and balancing workloads by reducing idle time and increasing the number of trips served. The adaptive power-of-$d$ approach turned out to be particularly effective, dynamically adapting to changes in demand and maintaining fleet performance over time. Ultimately, adaptive approaches demonstrated their ability to better manage system efficiency by adjusting vehicle availability and selection to yield the best overall performance across diverse scenarios.

\noindent The study also highlighted the importance of combining dynamic vehicle dispatch with efficient charging algorithms to ensure both high service levels and operational efficiency. Charging strategies that prioritize idle vehicles during off-peak hours and balance energy consumption throughout the day further improved service levels and reduced pickup times by optimizing fleet availability and charging frequency.

\noindent Overall, the findings suggest that achieving a balance between rapid response times and effective fleet management greatly enhances the performance of urban EV fleets. Using efficient matching and charging algorithms not only reduces idle times but also ensures that vehicles remain available to meet passenger demands without compromising energy management. The proposed framework offers a solid foundation for further research into optimizing EV operations in busy urban settings.

\subsection{Future Work}

\noindent Future research can build on the above findings by exploring more sophisticated matching and charging algorithms, as well as integrating additional data sources and real-world factors. Key areas for further study include:

\noindent \textbf{Additional Matching and Charging Algorithms:} Future simulations could investigate the impact of batching algorithms, where multiple requests are matched at once to reduce idle times and minimize fleet congestion. Adaptive batching strategies could be explored to dynamically group trip requests based on proximity and time windows to further reduce wait times and balance vehicle usage.

\noindent \textbf{Incorporating Traffic Congestion Models:} Integrating traffic congestion considerations into the simulation is vital, as the routing decisions send multiple vehicles in the same direction, and this can result in traffic congestion, which ultimately affects pickup times and efficiency. Future work should account for congestion effects in the matching algorithms to improve overall efficiency and avoid delays during high-demand periods.

\noindent \textbf{Integrating More Datasets:} Incorporating open-source datasets from other cities, such as the Ride Austin dataset, would enable broader validation of the proposed methods across different geographies and fleet configurations. This would also facilitate a more generalized model that is applicable to various urban environments to refine decision-making.

\noindent \textbf{Dynamic Pricing Mechanisms:} Implementing dynamic pricing strategies for EV rides and charging services could be an important area of exploration. By testing various pricing models, such as surge pricing during peak times or discounts for off-peak charging, future research could evaluate their effects on fleet demand, operational efficiency, and charging behaviors.

\noindent \textbf{Advanced Vehicle and Charger Optimization:} Further improvements can be made by refining the placement strategies for charging stations by utilizing predictive analytics to optimize fleet sizes and charger distribution in real time based on forecasted demand patterns.

\noindent By addressing these areas, further research can enhance the capabilities of the simulation framework, enabling it to tackle the growing challenges of urban EV fleet management. Advanced algorithms, real-world data integration, and additional operational constraints will ensure that the simulation framework evolves to support the increasing demand for electric vehicle services in dynamic urban settings.